%% file: Main_paper_deterministic.tex
\documentclass[smallextended,referee,envcountsect,]{svjour3}

\smartqed
\usepackage{graphicx}
\usepackage{amssymb,amsfonts, amsmath} %
 \usepackage{booktabs,amsmath,array} 
\usepackage{bm}
\usepackage{todonotes}
\usepackage{mathtools, relsize}
\usepackage{mathrsfs}
\usepackage{multirow,booktabs}
\usepackage{cancel}
\usepackage{xcolor}
\usepackage{url}
\usepackage{pgfplots}
\usepackage{tikz}
\usepackage{cleveref}
\usetikzlibrary{pgfplots.groupplots}
\pgfplotsset{compat=1.13}

\makeatletter

\newcommand{\R}{\ensuremath{\mathbb{R}}}

\newcommand{\rank}{\text{rank}}

\newcommand{\ML}{M}

\usepackage{enumitem}
\setlist[itemize]{topsep=-2pt,itemsep=0pt,leftmargin=15.0pt}

\usepackage{pgfplots}
\usepackage{graphicx,epstopdf} 
\usepackage{tikz}
\usepackage{subcaption}
\usepackage{algorithm}
\usepackage{algpseudocode}
\usepackage{algorithmicx}
\usepackage{caption}
\usepackage[toc,page]{appendix}

\usepackage{xspace}

\newcommand{\frob}{{\sc{frob}}\xspace}
\newcommand{\frobbnd}{{\sc{frobbnd}}\xspace}
\newcommand{\frobbndrk}{{\sc{frobbnd-$r_k$}}\xspace}

\newcommand{\orig}{{\sc{default}}\xspace}
\newcommand{\origtol}{{\sc{default-relaxed}}\xspace}

\usetikzlibrary{shapes.geometric}

\definecolor{matlabred}{rgb}{0.9047,    0.1918,    0.1988}
\definecolor{matlabblue}{rgb}{0.2941    0.5447    0.7494}
\definecolor{matlabgreen}{rgb}{	0.3718    0.7176    0.3612}
\definecolor{matlaborange}{rgb}{1.0000    0.5482    0.1000}

\usepackage{algpseudocode}
\algnewcommand{\LineComment}[1]{\Statex \(\%\) \small \textit{#1} \(\%\)}

\Crefname{ALC@unique}{Line}{Lines}
\usepackage{algorithm}
\DeclareMathOperator{\argmin}{arg\,min}

\begin{document}
\makeatletter
\renewcommand\subsection{\@startsection{subsection}{2}{\z@}%
  {-18pt \@plus -4pt \@minus -4pt}
  {6pt \@plus 2pt \@minus 2pt}
  {\normalfont\normalsize\bfseries\raggedright}}
\makeatother

\title{Convergence Analysis of the Inexact Singular Value Thresholding Algorithm}

\author{Stefania Bellavia \and
        Lisa Bertini \and \\
        Davide Palitta \and
        Margherita Porcelli
}

\institute{Stefania Bellavia \at
              Dipartimento di Ingegneria Industriale, Universit\`a degli Studi di Firenze, \\
              Viale Morgagni, 40/44, 50134 Florence, Italy \\
              \email{stefania.bellavia@unifi.it}\\
              Member of the INdAM Research Group GNCS and of the SIMAI-OPTIMA Activity Group.
           \and
           Lisa Bertini \at
              Dipartimento di Ingegneria dell'Informazione, Universit\`a degli Studi di Firenze, \\
              Via di S. Marta 3, 50139 Florence, Italy \\
              \email{lisa.bertini@unifi.it}\\
              Member of the INdAM Research Group GNCS.
           \and
           Davide Palitta \at
              Dipartimento di Matematica and (AM)$^2$, Alma Mater Studiorum Universit\`a di Bologna, \\
              Piazza di Porta San Donato 5, I-40127 Bologna, Italy \\
              \email{davide.palitta@unibo.it}\\
              Member of the INdAM Research Group GNCS.
           \and
           Margherita Porcelli \at
           Dipartimento di Ingegneria Industriale, Universit\`a degli Studi di Firenze, \\
              Viale Morgagni, 40/44, 50134 Florence, Italy, \\
              ISTI--CNR, Via Moruzzi 1, Pisa, Italy \\
              \email{margherita.porcelli@unifi.it}\\
              Member of the INdAM Research Group GNCS and of the SIMAI-OPTIMA Activity Group.
}


\maketitle

\begin{abstract}
The Singular Value Thresholding algorithm (SVT) is one of the most well-established methods for matrix completion, namely the problem of reconstructing a full matrix by observing only a handful of its entries. In spite of its easy implementation, SVT faces a serious bottleneck: the computation of an SVD at each iteration.
To lower this computational cost, existing SVT implementations often adopt iterative methods for computing an SVD-like factorization of the current iterate. Clearly, this strategy is able to compute only an approximation to the (truncated) SVD of the original matrix. To the best of our knowledge, this inexactness has never been accounted for in the converge analysis of SVT and the theoretical results available in the literature show 
the convergence of SVT only by assuming the SVDs to be computed exactly. This paper aims at filling this gap. In particular, we provide the right theoretical framework to show the convergence of SVT in case of inexact SVDs. Moreover, we specialize our results to the case where the SVDs are approximated by the Lanczos method and provide numerical results supporting our theoretical findings.

\vspace{0.2in} 

\end{abstract}



\section{Introduction}
In this work, we consider the solution of the so-called {\em matrix completion} problem. Roughly speaking, this amounts to the problem of reconstructing
 a matrix $A \in \mathbb{R}^{n_1 \times n_2}$ having at our disposal only a subset of its entries. Recovering a matrix from only $m$ observed entries is inherently impossible without additional information or structural assumptions, especially if $m \ll n_1 n_2$. In many applications, however, the target matrix $A$ is known to be low-rank, i.e., $\rank(A) \ll \min\{n_1, n_2\}$. Under this assumption, the matrix completion problem can be framed as follows:
\begin{equation}\label{matcompl-rank}
    \begin{aligned}
    \min_{X\in\mathbb{R}^{n_1\times n_2}} \quad \text{rank}(X) \quad \text{s.t.} \quad \mathcal{P}_{\Omega}(X) = \mathcal{P}_{\Omega}(A),
    \end{aligned}
\end{equation}
where $P_\Omega(\cdot)$ denotes the orthogonal projector onto the space of matrices that are zero outside the set $\Omega$, where, in our case, $\Omega$ denotes the set of indices of the entries of $A$ that can be observed,  namely $\Omega = \{ (i,j) \ | \ A_{ij} \mbox{ known} \}$.

Matrix completion problems are pervasive across various real-world applications, spanning from recommendation systems \cite{survey1} to time series imputation \cite{silei2023recovering}, and numerous approaches have been developed to address this challenge. A noncomplete list includes Riemannian optimization methods \cite{boumal2015low,goyens2023nonlinear,vandereycken2013low}, methods based on the formulation of (\ref{matcompl-rank}) as a nonconvex feasibility problem \cite{bellavia2025inexact,wang2021fast}, strategies based on matrix factorizations \cite{Relu26,kim2020deep,LMaFit}, and methods employing the reformulation as a  nuclear norm minimization problem \cite{recht2010guaranteed}; see also the survey papers \cite{survey1,survey2}
and the references therein.

The reformulation of~\eqref{matcompl-rank} based on nuclear norm minimization 
sees first the replacement of the objective function $\text{rank}(X)$ in~\eqref{matcompl-rank} with the nuclear norm of $X$, that is, solving 
the following relaxed optimization problem
\begin{equation}\label{matcompl-nn}
 \begin{aligned}
 \min_{X\in\mathbb{R}^{n_1\times n_2}} \quad \| X \|_*  \quad \text{s.t.} \quad \mathcal{P}_{\Omega}(X) = \mathcal{P}_{\Omega}(A),
   \end{aligned}
\end{equation} 
in place of~\eqref{matcompl-rank}, where $\|X \|_*:=\sum_{i=1}^{\text{rank}(X)}\sigma_i(X)$ and $\sigma_i(X)$ denotes the $i$-th singular value of $X$; see, e.g.,~\cite{CandesRecht2009}.
Notice that~\eqref{matcompl-nn} consists of the convex relaxation of the original problem~\eqref{matcompl-rank}. Cand\`es and Recht proved in \cite{CandesRecht2009} that if $\Omega$ 
is sampled uniformly at random among all subset of cardinality $m$ 
then, with high probability, the unique solution to (\ref{matcompl-nn}) 
is exactly $A$, provided that the number of samples obeys 
$m \ge \Gamma n^{5/4} \text{rank}(A) \log  n$, for some constant $\Gamma>0$ and $n:=\max\{n_1,n_2\}$.
In other words, under these assumptions, the problem in (\ref{matcompl-nn}) is ``formally equivalent'' 
to (\ref{matcompl-rank}).

A number of techniques based on the formulation (\ref{matcompl-nn}) can be found in the literature; see, e.g., gradient-type algorithms
\cite{cai2010singular,fpca,apgl}
and specialized second order methods for Semidefinite Programming formulation of (\ref{matcompl-nn}) \cite{jacek1,jacek2,loraine,liu2010interior,wang2025solving}.

In this work, we focus on the \textit{Singular Value Thresholding} (SVT) algorithm~\cite{cai2010singular}, which is one of the most prominent and well-established methods for the solution of~\eqref{matcompl-nn}. While this method has been extensively studied in the last decade, we believe there is still 
room for improving its theoretical understanding. 
In particular, the most demanding step of SVT is the computation of an SVD at each iteration.
Since computing a full SVD is often prohibitive\footnote{The computation of the full SVD of a matrix $Y\in\mathbb{R}^{n_1\times n_2}$ costs $\mathcal{O}(n_1n_2\min\{n_1,n_2\})$ floating point operations.}, the adoption of iterative methods is the only viable option in practice, giving rise to an inexact variant of the SVT algorithm. 
While the use of the latter can be traced back to the original SVT paper~\cite{cai2010singular}, available theoretical results are not able to guarantee the convergence of SVT under the computation of possibly inexact iterates. 
The main contribution of this work fills this gap by framing the right set of assumptions ensuring a convergent inexact SVT algorithm. 

The SVDs needed within SVT are often approximated by Krylov subspace methods; Golub–Kahan–Lanczos bidiagonalization~\cite{GolubKahan65} (abbreviated as
Lanczos bidiagonalization throughout the rest of this paper) is one of the most prominent methods in this class of numerical schemes. We thus specialize the general analysis mentioned above to the case of Lanczos bidiagonalization. In particular, we devise a couple of stopping criteria to equip Lanczos with.
A panel of numerical examples illustrates how these stopping criteria are such that the conditions of our analysis are met along with their computational burden.

Throughout the paper we adopt the following notation. Lower case letter ($x$) will denote vectors, upper case letters ($X$) matrices, and lower case Greek letters ($\alpha$) scalars.

For $x \in \mathbb{R}^p$, $\operatorname{diag}(x)$ denotes the diagonal matrix whose main diagonal contains the entries of $x$. The matrix inner product is defined as 
$\langle A, B \rangle = \operatorname{trace}(B^T A)$ for matrices of conforming dimensions and $\|X\|_F^2:=\langle X, X \rangle$ is the Frobenius norm of $X$. $\|X\|_2$ denotes the spectral norm of $X$ and $\sigma_i(X)$ is the $i$th singular values of the latter matrix.

\subsection*{Contributions}
The contributions of this work can be summarized as follows:
\vspace*{0.1cm}

\begin{itemize}
\item General convergence analysis of SVT with inexact SVD computations. 
  \item Specialization of the analysis to the case where the Lanczos method is adopted to approximate the SVDs.
        \end{itemize}
\vspace{0.3cm}
\noindent
Here is a synopsis of the paper. In section~\ref{sec_preliminaries} we report some background material. In particular, we revisit the standard SVT algorithm and its properties in section~\ref{sec_SVT}. Similarly, in section~\ref{sec_Lanczos} we recall the Lanczos method for SVD approximation. Section~\ref{sec_inexactSVT} sees the main contributions of this paper. Indeed, in this section we describe the right theoretical framework to show the convergence of SVT subject to the computation of inexact SVDs with generic error matrices (section~\ref{sec:deterministic_error}) whereas in section~\ref{sec_SVDwithLanczos} we specialize our analysis to the common case where the Lanczos method is adopted as iterative procedure to compute the approximate SVDs.
Some numerical results illustrating our findings are given in section~\ref{sec_Numerical_results} and the paper ends with section~\ref{sec_Conclusions} where we draw some conclusions. 

\section{Preliminaries}\label{sec_preliminaries}
In this section, we review the foundational principles of the SVT algorithm and of the Lanczos bidiagonalization process for SVD computation. These concepts are instrumental for introducing providing the necessary framework for the convergence study developed in this work.

\subsection{The SVT algorithm}\label{sec_SVT}
In this section, we outline the SVT algorithm for the matrix completion problem (\ref{matcompl-nn}), as originally proposed in \cite{cai2010singular}. We first introduce the fundamental shrinkage operator, which serves as the core building block of this iterative scheme.

 \begin{definition}
 Given a matrix $A \in \mathbb{R}^{n_1 \times n_2}$, let $\sigma_i(A)$ be the $i$th singular value of $A$. Given a parameter $\tau > 0$, the \textit{singular value thresholding operator} (or shrinkage operator) $\mathcal S_\tau:\mathbb{R}^{n_1 \times n_2}\rightarrow\mathbb{R}^{n_1 \times n_2}$ is defined as
 \begin{equation}\label{Operator_shrinkage}
    \mathcal S_\tau(A) = U\mathcal S_\tau(\Sigma) V^T, \qquad
     \mathcal S_\tau(\Sigma) := \operatorname{diag}((\sigma_1(A) - \tau)_+,\ldots,(\sigma_n(A) - \tau)_+),
 \end{equation}
 where $n:=\min\{n_1,n_2\}$ and $(t)_+$ denotes $\max(t,0)$.  
 
\end{definition}

For $\tau>0$, the SVT algorithm seeks the unique solution to the following closely related problem
\begin{equation}\label{PB_SVT}
\begin{aligned}
    \min_{X\in\mathbb{R}^{n_1\times n_2}}& \ \ f_\tau(X)=\tau\|X\|_*+\frac{1}{2}\|X\|^2_F  \quad \text{s.t.} \quad \mathcal{P}_{\Omega}(X) = \mathcal{P}_{\Omega}(A).
\end{aligned}
\end{equation}
The solution of \eqref{PB_SVT} converges to the solution of the formulation \eqref{matcompl-nn} with minimum Frobenius norm as $\tau$ tends to infinity, as shown in 
\cite[Theorem 3.1]{cai2010singular}.

The SVT method is an iterative algorithm that, given $\tau>0$, a sequence of positive step sizes $\{\delta_k\}_{k\geq1}$, and an initial guess $Y^0$, performs the following steps at each iteration $k$
\begin{equation}\label{Sequences_SVT_alg}
\begin{cases}
X^k = \mathcal S_\tau(Y^{k-1}), \\
Y^k = Y^{k-1} + \delta_k \bigg(\mathcal P_\Omega(A) - \mathcal P_\Omega(X^k)\bigg).
\end{cases}
\end{equation}
The iterates $\{X^k\}$ and $\{Y^k\}$ defined in  \eqref{Sequences_SVT_alg} 
satisfy two important properties.
\begin{itemize}
    \item \textbf{Low rank of $X^k$}. Although there is no theoretical proof in the literature, the matrices $\{X^k\}$ often show a very low rank. This phenomenon is easy to understand: by using a large threshold $\tau$, the $\mathcal S_\tau(\cdot)$ operator sets most singular values of $X^k$ to zero, and therefore $X^k=\mathcal S_\tau(Y^k)$ has low rank. This is a significant feature, given that the target solution $A$ is expected to be low-rank. 
    \item\textbf{Sparsity of $Y^k$}. Provided $\mathcal{P}_\Omega(Y^0)=Y^0$, the matrices $\{Y^k\}$ all have the sparsity pattern encoded in $\Omega$. The smaller number $m$ of available entries, the higher the degree of sparsity in the matrix $Y^k$. This property is crucial, as sparsity can be exploited to optimize memory storage and significantly accelerate computations.
\end{itemize}
\vskip 0.3cm
\noindent
The converge property of the SVT method is recalled in the theorem below, for the proof we refer the reader to~\cite[Theorem 4.2]{cai2010singular}.

\begin{theorem}\label{thm:Conv_SVT}
If the step sizes $\delta_k$ satisfies the condition
\begin{equation}\label{eq:deltak}
0 < \inf \delta_k \leq \sup \delta_k < 2,    
\end{equation}
then the sequence $\{X^k\}$ obtained through \eqref{Sequences_SVT_alg} converges to the unique solution of \eqref{PB_SVT}.
\end{theorem}

We report in Algorithm \ref{alg:svt} the details of the SVT algorithm 
whereas the application of the soft thresholding operator is depicted in Algorithm \ref{alg:StauSVT}.

\begin{algorithm}[t!]
    \caption{Singular Value Thresholding Algorithm (SVT)}
    \label{alg:svt}
    \begin{algorithmic}[1]
        \Require Sample set $\Omega$ and observed entries $\mathcal P_{\Omega}(A)$, a sequence of step sizes $\{ \delta_k \}$, threshold parameter $\tau$, tolerance $tol$, max. iteration count $k_{max}$, initial guess $Y^0$
        \Ensure $X^{opt}$
        \For {$k = 1, 2, \dots, k_{max}$}
                \State Compute $\mathcal S_{\tau}(Y^{k-1})=U^{k-1}\Sigma^{k-1}(V^{k-1})^T$
                  \Comment{via Algorithm \ref{alg:StauSVT} applied to }\label{alg:subroutine}
                \Statex      \Comment{$Y^{k-1}$ and threshold $\tau$}
            \State Set $r_k := \max\{j : \sigma_j\left(\mathcal S_\tau(Y^{k-1})\right) > 0\}$
            \State Set $X^k := U^{k-1}(:,1:r_k)\Sigma^{k-1}(1:r_k,1:r_k)(V^{k-1}(:,1:r_k))^T$
            \If {$\|\mathcal P_\Omega(X^k) - \mathcal P_\Omega(A)\|_F / \|\mathcal P_\Omega(A)\|_F < tol$} 
                \State \textbf{break}
            \EndIf
            \State Set $Y^k = Y^{k-1} + \delta_k \bigg(\mathcal P_\Omega(A) -\mathcal P_\Omega(X^k)\bigg)$
        \EndFor
        \State \textbf{return} $X^{opt} = X^k$
    \end{algorithmic}
\end{algorithm}

\begin{algorithm}[t!]
    \caption{Soft-thresholding operator evaluation with exact SVD}
    \label{alg:StauSVT}
    \begin{algorithmic}[1]
        \Require Matrix $Y$, threshold $\tau$
        \Ensure the SVD decomposition of $\mathcal S_\tau(Y)$ 
       \State Compute the SVD $Y = U\Sigma V^T$
       \State Set $\Sigma=(\Sigma-\tau I)_+$
       \State \textbf{return} $[U, \Sigma, V]$
    \end{algorithmic}
\end{algorithm}

From a computational perspective, evaluating the full SVD at each iteration of Algorithm \ref{alg:svt} becomes prohibitively expensive, or even unfeasible, when dealing with large-scale matrices. To overcome this limitation, the full SVD can be replaced by an iterative, truncated SVD procedure, which computes only the first largest singular values and their corresponding singular vectors. However, since the effective rank of the solution is not known a priori, fixing a static number of singular values to compute might lead to inaccurate thresholding if this number is chosen to be smaller than the rank of the sought solution.

To address this issue, an adaptive strategy, originally introduced in \cite{cai2010singular}, is implemented in practice; see section \ref{sec_Numerical_results}. 

\subsection{Lanczos bidiagonalization for SVD approximations}\label{sec_Lanczos}
 Given a matrix $Y\in\R^{n_1\times n_2}$ (this will be our SVT iterate $Y^{k-1}$) Lanczos bidiagonalization iteratively constructs an approximation to its SVD.
In particular, by selecting a vector $p_1\in\mathbb{R}^{n_1}$ with unit norm, the zero vector $q_0\in\mathbb{R}^{n_2}$, and $\beta_1=1$, the following recurrence relations define Lanczos bidiagonalization
\begin{equation}\label{eq:rel_lanczos1}
      \alpha_\ell q_\ell =Y^T p_\ell - \beta_\ell q_{\ell-1}, 
    \qquad \beta_{\ell+1} p_{\ell+1} = Y q_\ell - \alpha_\ell p_\ell,\quad \ell\geq1, %
 \end{equation}    
where $\alpha_\ell, \ \beta_\ell$ are chosen such that the vectors $p_\ell$'s and $q_\ell$'s, called \textit{Lanczos vectors}, have unitary norm.
If we consider the matrices {$P_{\ell} = [p_1, \dots, p_{\ell}]\in\mathbb{R}^{n_1\times \ell}$} and $Q_\ell = [q_1, \dots, q_{\ell}]\in\mathbb{R}^{n_2\times \ell}$ it can be shown that these have orthonormal columns, namely $P_\ell^TP_\ell=I$ and $Q_\ell^TQ_\ell=I$, and their column spaces are indeed Krylov subspaces. In particular, $\text{Range}(Q_\ell)=\mathcal{K}_\ell(Y^TY,Y^Tp_1)$ and $\text{Range}(P_\ell)=\mathcal{K}_\ell(YY^T,p_1)$. Moreover, if we collect the coefficients $\alpha_\ell$'s and $\beta_\ell$'s in the following bidiagonal matrix
\begin{equation}\label{eq:Bell_Lanczos}
    B_\ell = 
    \left(
    \begin{array}{c@{\hspace{.5em}}c@{\hspace{.5em}}c@{\hspace{.5em}}c}
        \alpha_1 &           &        &  \\
        \beta_2  & \alpha_2  &        &  \\
                 & \beta_3   & \smash{\ddots} &  \\
                 &           & \smash{\ddots} & \alpha_\ell \\
                 &           &        & \beta_{\ell+1}
    \end{array}
    \right) \in \R^{(\ell+1) \times \ell},
\end{equation}
then we can rewrite the recurrences~\eqref{eq:rel_lanczos1} in matrix form as
\begin{equation}\label{eq:rel_lanczos2}
    Y Q_\ell = P_{\ell + 1} B_\ell, \qquad  Y^T P_{\ell +1} = Q_\ell B_\ell^T + \alpha_{\ell+1}q_{\ell+1}e_{\ell+1}^T, 
\end{equation}
where $e_{\ell+1}$ is the $(\ell+1$)-th vector of the canonical basis. Recalling the orthogonality of the columns of $P_{\ell+1}$, we thus have $B_\ell=P_{\ell+1}^TYQ_\ell$ so that the matrix $ \ML_\ell= {P_{\ell+1} B_{\ell} Q_{\ell}^T}=({P_{\ell+1} P_{\ell+1}^T)Y(Q_\ell Q_{\ell}^T)}\in\R^{n_1\times n_2}$ may amount to a good approximation to the original $Y$. In particular, we are interested in how the SVD of $\ML_\ell$ gives a good approximation to that of $Y$.
First of all, the SVD of $\ML_\ell$ is easy to construct. Indeed, if $B_\ell = U_B \Sigma_B V_B^T$ is the SVD of the small dimensional $B_\ell$, then
\begin{equation} \label{svd_Z}
    \ML_\ell = {P_{\ell+1} B_{\ell} Q_{\ell}^T}= \underbrace{P_{\ell+1} U_B}_{\tilde U_\ell} \underbrace{\Sigma_B}_{\tilde \Sigma_\ell} \underbrace{V_B^T Q_\ell^T}_{\tilde V_\ell^T}=\tilde U_\ell \tilde \Sigma_\ell \tilde V_\ell^T,
\end{equation}
gives the SVD of $\ML_\ell$, where $\tilde \Sigma_\ell{\in\mathbb{R}^{(\ell+1)\times \ell}}$ is a rectangular diagonal matrix whose diagonal entries $\tilde \sigma_{1,\ell}(M_\ell), \dots, \tilde \sigma_{\ell,\ell}(M_\ell)$ can be seen as approximations of the largest $\ell$ singular values {$\sigma_1(Y),\dots,\sigma_\ell$(Y)} of $Y$. 
 Notice that, in case  there exists $\hat j \in \{{1},\ldots,\ell\}$ such  that $\alpha_{\hat j}=0$, Lanczos bidiagonalization breaks down and  $\tilde{\sigma}_{i,\hat j}(M_\ell)$,  $i=1,\ldots,\hat j$ are  exact singular values of $ Y$   and  $\tilde V_{\hat j}$ is an invariant subspace \cite[Chapter IV]{saad2011numerical}.
{If $\beta_{\ell+1} = 0$, then from the equations in~ \eqref{eq:rel_lanczos2} it follows that
$$
 Y^T P_{\ell} = Q_{\ell} B_\ell^\top(1:\ell,1:\ell),\qquad
	Y Q_\ell = P_\ell B_\ell(1:\ell,1:\ell),
$$
where {$B_\ell(1:\ell,1:\ell)$ denotes the principal square $\ell\times\ell$ submatrix of $B_\ell$. 
Therefore,
\[
Y^{T} Y Q_\ell =  Y^{T} P_\ell B_\ell(1:\ell,1:\ell) = Q_{\ell} B_\ell^T(1:\ell,1:\ell) B_\ell(1:\ell,1:\ell),
\]
which means that the singular values of  $B_\ell(1:\ell, 1:\ell)$ are exact singular values of $ Y$.}}

The preceding discussion indicates how singular values estimates can be obtained via the Lanczos method, but it reveals nothing about the approximation quality of the singular values and the corresponding singular vectors. 
We refer to \cite[Section 10.4]{golub2013matrix} and references therein for the convergence results of the Lanczos procedure. Here, we highlight the following properties that are crucial for our analysis,  {in the common situation where the algorithms does not break down for null $\alpha_j$ or $\beta_j$}: 
\begin{itemize}
    \item in exact arithmetic the method  has finite termination; for $\ell =\min\{n_1,n_2\}$, the SVD in \eqref{svd_Z} is the exact SVD of $Y$. However, the procedure is generally implemented as an iterative procedure due to round-off errors and to avoid to perform a large number of steps;
    
    \item the inexact singular values $\tilde \sigma_{i,\ell}(M_\ell)$ approximate the corresponding $\sigma_i(Y)$ from below \cite[Theorem  6.4]{saad2011numerical};
    \item letting $(U_B)_{\ell +1, j}$ be the last element of the $j$-th column of $U_B$ in \eqref{svd_Z}, as shown in \cite{propack_doc}, we have that 
    \begin{equation} \label{eq:Lanczos_error_bound_v2}
        \min_{\sigma(Y) \in \Sigma (Y)} | \sigma(Y) - \tilde \sigma_{j, \ell}(M_\ell) | \le | \alpha_{\ell+1} | |(U_B)_{\ell +1, j}|.
    \end{equation}
{\item the singular values $\tilde \sigma_{1,\ell}(M_\ell), \dots, \tilde \sigma_{\ell,\ell}(M_\ell)$ of $M_\ell$
    computed by the Lanczos method are distinct, independently on the singular values distribution of $Y$, as they are the singular values of the bidiagonal matrix  $B_\ell$ with nonzero entries; 
        see \cite[Theorem 2.5]{higham_power_2024}.}
\end{itemize}
\vspace{.3cm}

The latter bullet says that, in practice, the Lanczos method acts as if all singular values of $Y$ were simple \cite{li_convergence_2015}, even if this may not be the case. To approximate all or some of the copies of a multiple singular value, one possibility is to  use a \emph{block} Lanczos method with a block size greater than $j$, which requires prior knowledge of the distribution of the singular values of $\widetilde Y^{k-1}$ \cite{li_convergence_2015}. While this is clearly possible, to simplify our derivations, in the following we will assume all the singular values of $Y$ to be distinct as done in~\cite{saad_rates_1980}.

\section{The inexact SVT method}\label{sec_inexactSVT}
 In this section, we analyze the convergence of the SVT method when the computation of the soft-thresholding operator is affected by error yielding in fact an inexact version of SVT. 

\subsection{Convergence analysis}\label{sec:deterministic_error}
We consider a scenario where the SVT sequences have the form 
\begin{equation}\label{Successioni_SVT_calc}
\begin{cases}
    \widetilde{X}^k&=\mathcal S_\tau(\widetilde{Y}^{k-1})+\Delta^k,\\
    \widetilde{Y}^k&=\widetilde{Y}^{k-1}+\delta_k\bigg(\mathcal{P}_\Omega(A)-\mathcal{P}_\Omega(\widetilde{X}^k)\bigg),
\end{cases}
\end{equation}
where we denote by $\Delta^k$ the error introduced in the shrinkage operator calculation when computing the SVD. Notice that in this case the sequence $\widetilde{Y}^k$ is defined as in \eqref{Sequences_SVT_alg}. The theorem below  generalizes Theorem \ref{thm:Conv_SVT} and provides a bound on the limit of the averaged error. 

\begin{theorem}
\label{thm:Conv_SVT_errore_det}
  Let $\{\widetilde{X}^k\}$ be the sequence defined in \eqref{Successioni_SVT_calc} and $X^*$ be the solution of the problem \eqref{PB_SVT}. If the sequence $\{\delta_k\}$ satisfies (\ref{eq:deltak}) and there exists $\theta_*>0$ such that 
  $\|\Delta^k\|_*\leq\theta_*$ for any $k=0,1,\ldots$, then
\begin{equation} \label{eq:det_decrease}
   \lim_{K\to+\infty} \sum_{k=0}^K\frac{\|\widetilde{X}^k-X^*\|_F^2}{K+1} \leq \Bigg(\frac{64\gamma^2}{\beta^2}+1\Bigg)\theta_*^2
+\frac{32\tau\gamma}{\beta}\theta_*, 
\end{equation}
where $\beta=\inf_k(\,2\delta_k-\delta_k^2)>0$ and $\gamma=\sup_k\delta_k$.

\end{theorem}
\begin{proof}
Consider the Lagrangian associated with problem \eqref{PB_SVT}
\begin{equation}\label{Lagrangiana}
    \mathcal{L}(X,Y)=f_{\tau}(X)+\langle Y, P_{\Omega}(A-X) \rangle,
\end{equation}
where $Y$ is the Lagrange multiplier associated with the equality constraint.
Let $X^*$ be the solution of the same constrained minimization problem~\eqref{PB_SVT} and $Y^*$ be the corresponding Lagrangian multiplier. Then, the couple  $(X^*,Y^*)$
 clearly satisfies the necessary optimality condition of the problem, which includes being a stationary point for the associated Lagrangian defined in (\ref{Lagrangiana}), that is $\nabla_X\mathcal{L}(X^*,Y^*)=0$ which, by explicit computation, is equivalent to 
\begin{equation}\label{eq:Z*}
    0=Z^*-\mathcal P_\Omega(Y^*), \quad \text{with } Z^*\in \partial f_\tau(X^*),
\end{equation}
where $ \partial f_\tau$ denotes the subdifferential of the function $f_\tau$.
Moreover, by construction of the sequence \eqref{Successioni_SVT_calc} and considering that for fixed $Y$, the minimizer of the Lagrangian $\mathcal{L}(X,Y)$ is given by the singular value shrinkage of $P_\Omega(Y)$, i.e. $\operatorname{argmin}_X \mathcal{L}(X,Y) = \mathcal S_\tau(P_\Omega(Y))$ \cite[Sec. 2.4]{cai2010singular}, we have that for every $k=1,2,\dots$, $\widetilde{X}^k-\Delta^k$ minimizes the Lagrangian with $\widetilde{Y}^{k-1}$ fixed; thus this matrix also satisfies the same condition, that is $\nabla_X\mathcal{L}(\widetilde X^k-\Delta^k,\widetilde Y^{k-1})=0$, which, by explicitly evaluating the Lagrangian gradient, yields exactly
\begin{equation}\label{eq:Ztilde}
0=\widetilde{Z}^k-\mathcal P_\Omega(\widetilde{Y}^{k-1}), \quad \text{with } \widetilde{Z}^k\in \partial f_\tau(\widetilde{X}^k-\Delta^k).
\end{equation}
Subtracting the equation (\ref{eq:Z*}) from equation (\ref{eq:Ztilde}) yields
$    (\widetilde{Z}^k-Z^*)-\mathcal P_\Omega(\widetilde{Y}^{k-1}-Y^*)=0.
$

By \cite[Lemma 4.1]{cai2010singular}, we obtain
\begin{align} \label{bound_scalar_prod_short}
\langle \widetilde{X}^k-X^*, \mathcal P_\Omega(\widetilde{Y}^{k-1}-Y^*)\rangle&= \langle \widetilde{X}^k-X^*, \widetilde{Z}^k-Z^*\rangle\nonumber\\ 
&\geq \|\widetilde{X}^k-\Delta^k-X^*\|_F^2+\langle \widetilde{Z}^k-Z^*,\Delta^k\rangle.
\end{align}
Since $(X^*,Y^*)$ is the optimal solution of \eqref{PB_SVT}, clearly $X^*$ satisfies the constraint $\mathcal P_\Omega(A)=\mathcal P_\Omega(X^*)$. Therefore, using the definition of the sequence $\{\widetilde{Y}^k\}_{k\geq 0}$ in \eqref{Successioni_SVT_calc}, we have
$\|\mathcal P_\Omega(\widetilde{Y}^k-Y^*)\|_F
=\|\mathcal P_\Omega(\widetilde{Y}^{k-1}-Y^*)+\delta_k\mathcal P_\Omega(X^*-\widetilde{X}^k)\|_F.
$
Define $\widetilde{r}_k=\|\mathcal P_\Omega(\widetilde{Y}^k-Y^*)\|_F$. By using the inequality \eqref{bound_scalar_prod_short}, we obtain
\begin{align}
\widetilde{r}_k^2
&=\widetilde{r}_{k-1}^2-2\delta_k\langle \mathcal P_\Omega(\widetilde{Y}^{k-1}-Y^*), \widetilde{X}^k-X^* \rangle+ \delta_k^2\|\mathcal P_\Omega(X^*-\widetilde{X}^k)\|_F^2 \notag \\
&\leq \widetilde{r}_{k-1}^2 -2\delta_k(\|\widetilde{X}^k-\Delta^k-X^*\|_F^2+\langle \widetilde{Z}^k-Z^*,\Delta^k\rangle) +\delta_k^2\|\widetilde{X}^k-X^*\|_F^2. \notag
\end{align}
since for every matrix $X$ we have $\|\mathcal P_\Omega(X)\|_F\leq\|X\|_F$. 

The term $\langle \widetilde{Z}^k-Z^*,\Delta^k\rangle$ is bounded using the Von Neumann's trace inequality \cite[Theorem~3.3.14]{horn1991topics}, the characterization of the subdifferential of the nuclear norm~\cite{cai2010singular}, Lemma \ref{stima_subdiff} and considering that $Z^*=\tau Z_0^*+X^*$ with $Z^*_0\in \partial\|X^*\|_*$ and similarly for $\widetilde{Z}^k$ (see Lemma \ref{lemma_char_Z}). We also use the Cauchy-Schwarz inequality
\begin{align}
|\langle \widetilde{Z}^k-Z^*,\Delta^k\rangle| &=|\tau\langle \widetilde{Z}_0^k-Z_0^*,\Delta^k\rangle+\langle \widetilde X^k-\Delta^k-X^*,\Delta^k\rangle|\nonumber\\
&\leq\tau\|\widetilde{Z}_0^k-Z_0^*\|_2\|\Delta^k\|_*+\| \widetilde X^k-\Delta^k-X^*\|_F\|\Delta^k\|_F\nonumber\\
&\leq 4\tau\|\Delta^k\|_*+\|\widetilde{X}^k-\Delta^k-X^*\|_F\|\Delta^k\|_F.\notag
\end{align}
For the term $\|\widetilde{X}^k-\Delta^k-X^*\|^2_F$, we apply the triangle inequality and get
\begin{align*}
\|\widetilde{X}^k-\Delta^k-X^*\|^2_F & \geq (\|\widetilde{X}^k-X^*\|_F-\|\Delta^k\|_F)^2 \nonumber\\
& = \|\widetilde{X}^k-X^*\|_F^2+\|\Delta^k\|_F^2 -2
\|\widetilde{X}^k-X^*\|_F\|\Delta^k\|_F.
\end{align*}
Overall we get 
\begin{align}
    \widetilde{r}_k^2&\leq \widetilde{r}_{k-1}^2
-(2\delta_k-\delta_k^2)||\widetilde{X}^k-X^*||_F^2
+4\delta_k||\widetilde{X}^k-X^*||_F||\Delta^k||_F \nonumber\\
&\quad +8\delta_k\tau||\Delta^k||_*+2\delta_k||\widetilde{X}^k-\Delta^k-X^*||_F||\Delta^k||_F-2\delta_k||\Delta^k||_F^2\notag\\
&\leq \widetilde{r}_{k-1}^2
-(2\delta_k-\delta_k^2)||\widetilde{X}^k-X^*||_F^2
+2\delta_k\left(\epsilon||\widetilde{X}^k-X^*||^2_F
+\frac{||\Delta^k||^2_F}{\epsilon}\right) \nonumber\\
&\qquad +8\delta_k\tau||\Delta^k||_*-2\delta_k||\Delta^k||_F^2
+\delta_k \left(\epsilon_1||\widetilde{X}^k-\Delta^k-X^*||^2_F
+\frac{||\Delta^k||^2_F}{\epsilon_1}\right).\notag
\end{align}
where in the second inequality we used  the Young's inequality\footnote{For any $ a,b$ and $\epsilon>0$, it holds $ab\leq\frac{1}{2}(\epsilon a^2+\frac{b^2}{\epsilon})$,} to the terms $||\widetilde{X}^k-X^*||_F||\Delta^k||_F$ and $||\widetilde{X}^k-\Delta^k-X^*||_F||\Delta^k||_F$ with $\epsilon$ and $\epsilon_1$,  respectively.

Moreover, letting $\beta = \inf_k(2\delta_k-\delta_k^2) > 0$, choosing the appropriate constants $\epsilon=\frac{\beta}{4\delta_k}$ and $\epsilon_1=\frac{\beta}{8\delta_k}$ and using the convexity of $\|\cdot\|_F^2$\footnote{For any matrix $A$ and $B$ we have $\|A+B\|_F^2\leq 2\|A\|_F^2+2\|B\|_F^2$.} for the term $||\widetilde{X}^k-\Delta^k-X^*||^2_F$, we obtain 
\begin{equation*}
\widetilde{r}_k^2 \leq \widetilde{r}_{k-1}^2 -\frac{\beta}{4}\|\widetilde{X}^k-X^*\|_F^2 +\left(\frac{16\delta_k^2}{\beta}+\frac{\beta}{4}\right)\|\Delta^k\|^2_F +8\delta_k\tau\|\Delta^k\|_*.
\end{equation*}
Summing $k$ from $0$ to $K$, noting that $\widetilde{r}_K^2 \geq 0$ and defining $\gamma=\sup_k \delta_k$, we get
\begin{equation}\label{eq:bound_deter_SVT}
\frac{\beta}{4}\sum_{k=0}^K\frac{\|\widetilde{X}^k-X^*\|_F^2}{K+1} \leq \frac{\widetilde{r}_0^2}{K+1} +\left(\frac{16\gamma^2}{\beta}+\frac{\beta}{4}\right)\sum_{k=0}^K\frac{\|\Delta^k\|^2_F}{K+1} +8\tau \gamma\sum_{k=0}^K \frac{\|\Delta^k\|_*}{K+1}.
\end{equation}
Since $\|\Delta^k\|_F\le \|\Delta^k\|_*$ and the  nuclear norm of the matrices $\Delta^k$ is uniformly bounded from above by $\theta_*$ by assumption,
as $K \to \infty$ the average error satisfies inequality \eqref{eq:det_decrease}.
\end{proof}

The inequality (\ref{eq:bound_deter_SVT}) in the proof of
 Theorem \ref{thm:Conv_SVT_errore_det} 
shows that, on average, the squared error  $\|\widetilde{X}^k-X^*\|_F^2$ is bounded by a term that goes to zero as $O(1/K)$ plus a constant term that depends on $\sum_{k=0}^K\|\Delta^k\|_F^2/(K+1)$ and $\sum_{k=0}^K\|\Delta^k\|_*^2/(K+1)$. Consequently, the right-hand side provides an explicit upper bound that acts as an asymptotic error floor, preventing the convergence of the infinite series of errors unless the series of errors are bounded.

\begin{corollary}\label{coll:svt_conv}
Under the hypotheses of Theorem \ref{thm:Conv_SVT_errore_det}, if we assume the stronger condition that 
\begin{equation} \label{strong_error_condition}
\sum_{k=0}^\infty \|\Delta^k\|^2_F < \infty \quad \text{and} \quad \sum_{k=0}^\infty \|\Delta^k\|_* < \infty,
\end{equation}
then the series $\sum_{k=0}^\infty \|\widetilde{X}^k-X^*\|_F^2$ converges and consequently $\lim_{k\to\infty} \|\widetilde{X}^k-X^*\|_F^2 = 0$.
\end{corollary}
\begin{proof}

The thesis follows from inequality \eqref{eq:bound_deter_SVT}.
Taking the limit as $K \to \infty$ on both sides we get
\begin{equation*}
\frac{\beta}{4} \sum_{k=0}^\infty \|\widetilde{X}^k-X^*\|_F^2 \leq \widetilde{r}_0^2 + \left(\frac{16\gamma^2}{\beta}+\frac{\beta}{4}\right) \sum_{k=0}^\infty \|\Delta^k\|^2_F + 8\gamma\tau \sum_{k=0}^\infty \|\Delta^k\|_*.
\end{equation*}
Since the series of the errors $\sum_{k=0}^\infty \|\Delta^k\|^2_F$ and $\sum_{k=0}^\infty \|\Delta^k\|_*$ are convergent and then bounded by assumption, it follows that the series $\sum_{k=0}^\infty \|\widetilde{X}^k-X^*\|_F^2$ converges. As a necessary condition for the convergence of a series, the general term must vanish, which implies $\lim_{k\to\infty} \|\widetilde{X}^k-X^*\|_F^2 = 0$.
\end{proof}
While Theorem \ref{thm:Conv_SVT_errore_det} only guarantees an asymptotic error floor under the assumption of bounded perturbations, with the stronger requirement of a convergent series of errors, Corollary \ref{coll:svt_conv} ensures full asymptotic convergence to the optimal solution. Formally, for the error series $\sum_{k=0}^\infty \|\Delta^k\|^2_F$ and $\sum_{k=0}^\infty \|\Delta^k\|_*$ to converge, the sequences of perturbations must vanish sufficiently fast (strictly faster than $\mathcal{O}(1/k)$). Under this assumption, the cumulative effect of the inexactness introduced in the SVD step remains bounded, allowing the error $\|\widetilde{X}^k-X^*\|_F^2$ to vanish asymptotically.

\subsection{The case of Lanczos bidiagonalization}\label{sec_SVDwithLanczos}

In this section we assume that the Lanczos bidiagonalization is employed within Algorithm~\ref{alg:StauSVT} to approximate the SVD of $\widetilde Y^{k-1}$. If $\ML_{\ell_k}$ denotes the matrix computed by Lanczos after $\ell_k$ iterations, then we can look at the error matrix $\Delta^k$ in~\eqref{Successioni_SVT_calc} as 
\begin{equation}\label{eq:defDeltak}
\Delta^k=\mathcal{S}_\tau(\ML_{\ell_k})-\mathcal{S}_\tau(\widetilde Y^{k-1}).    
\end{equation}

Our goal is to show how to choose the number of Lanczos iterations $\ell_k$ so that the conditions in~\eqref{strong_error_condition} are satisfied.
Since $\Delta^k$ depends on $\widetilde Y^{k-1}$ only through the truncation error of the Lanczos approximation $M_{\ell_k}$, we start by relating $\|\Delta^k\|_F$ to quantities that are available, or can be estimated, during the bidiagonalization.
By Lemma~\ref{est_op_shr_fr} and Lemma~\ref{error_bound_Ptau} we have
\begin{align}
 \|\Delta^k\|_F^2 & = \|\mathcal S_\tau(M_{\ell_k})-\mathcal S_\tau(\widetilde Y^{k-1})\|_F^2\nonumber  \\
 &\leq \|M_{\ell_k}-\widetilde Y^{k-1}\|_F^2-\|\mathcal P_\tau(M_{\ell_k})-\mathcal P_\tau(\widetilde Y^{k-1})\|_F^2\nonumber \\
    &=\|M_{\ell_k}-\widetilde Y^{k-1}\|_F^2-\sum_{i>\ell_k}\min(\sigma_i(\widetilde Y^{k-1}),\tau)^2.\label{eq:error_lanczos}
\end{align}
The right-hand side of~\eqref{eq:error_lanczos} is the quantity on which we build the stopping rule of the Lanczos process:
the conditions in~\eqref{strong_error_condition} are met as soon as it decays like $\frac{C}{k^{2+\epsilon}}$, for given constants $C,\epsilon>0$\footnote{Since $\|A\|_*^2\le\operatorname{rank}(A)\,\|A\|_F^2$ for any matrix $A$, the decay imposed on $\|\Delta^k\|_F^2$ also yield $\|\Delta^k\|_*\le\sqrt{\min(n_1,n_2)}\,\frac{\sqrt{C}}{k^{1+\epsilon/2}}$, so that both 
conditions in~\eqref{strong_error_condition} hold.}.

Starting from~\eqref{eq:error_lanczos}, we derive three stopping criteria. Two of them, \frob and \frobbnd, are rigorous upper bounds on $\|\Delta^k\|_F^2$, so that the convergence theory developed in the previous sections applies.
The last one, \frobbndrk is heuristic and exploits the fact that the soft-thresholding operator cuts off all the singular values below $\tau$.
In all cases, the Lanczos bidiagonalization is stopped at the first iteration $\ell_k$ at which the corresponding quantity falls below $\frac{C}{k^{2+\epsilon}}$ .

\paragraph{A conservative bound - \frob}
The recurrences of the Lanczos bidiagonalization can be employed to cheaply compute $\|M_{\ell_k}-\widetilde Y^{k-1}\|_F^2$ on the fly as shown in~\cite[Lemma 4.11]{bellavia2025inexact}. Indeed, if $\omega_{\ell_k}=\|M_{\ell_k}-\widetilde Y^{k-1}\|_F^2$, then it holds
\begin{equation}\label{eq:rec}
    \omega_{\ell_k-1}=\omega_{\ell_k}+\alpha_{\ell_k}^2+\beta_{\ell_k+1}^2,
\end{equation}
where $\alpha_{\ell_k}$ and $\beta_{\ell_k+1}$ are the quantities in~\eqref{eq:rel_lanczos1}. This means that 
 checking $\|M_{\ell_k}-\widetilde Y^{k-1}\|_F^2$ 
does not require any additional computational load. 
On the other hand, it is not as easy to derive lower bounds for the terms $\min(\sigma_i(\widetilde Y^{k-1}),\tau)^2$ in~\eqref{eq:error_lanczos} without further assumptions and we are not aware of any properties of the Lanczos algorithm able to provide computable estimates thereof, unless the whole space is generated, i.e., $\ell_k=\min(n_1,n_2)$.
However, we can consider the coarser bound 
$$\|\Delta^k\|_F^2 \leq \omega_{\ell_k},$$
and  run $\ell_k$ Lanczos iterations to approximate the SVD of $\widetilde Y^{k-1}$ such that 
\begin{equation}\label{eq:Lancsoz_stop}
    \omega_{\ell_k}<\frac{C}{k^{2+\epsilon}},
\end{equation}
for given constants $C,\epsilon>0$.

A couple of considerations about the Lanczos stopping criterion~\eqref{eq:Lancsoz_stop} are as follows. First, at the initial phase of SVT, we expect to run a few Lanczos iterations in contrast to what happens by adopting a fixed, restrictive tolerance.
While this certainly decreases the per-iteration cost of that initial SVT stage, it may also induce a delay in the overall SVT convergence. 
Second, solely looking at $\omega_{\ell_k}$ may still lead to running more Lanczos iterations than necessary. Indeed, in~\eqref{eq:error_lanczos} we completely disregarded the term $-\sum_{i>\ell_k}\min(\sigma_i(\widetilde Y^{k-1}),\tau)^2$.

\paragraph{Refined bound - \frobbnd}
The term $-\sum_{i>\ell_k}\min(\sigma_i(\widetilde Y^{k-1}),\tau)^2$ in~\eqref{eq:error_lanczos}  can be bounded by using the bound (\ref{eq:Lanczos_error_bound_v2}) if we assume that $\sigma_i(\widetilde Y^{k-1})\leq \tau$ for any $i>\ell_k$ and $\sigma_j(\widetilde Y^{k-1})=\argmin_{\sigma(\widetilde Y^{k-1})\in\Sigma(\widetilde Y^{k-1})}|\sigma(\widetilde Y^{k-1})-\tilde\sigma_{j,\ell_k}(M_{\ell_k})|$ for any $j\leq\ell_k$ where we remind the reader that $\sigma_j(M_{\ell_k})$ is the $j$th Ritz values computed at the $\ell_k$ Lanczos iteration. By denoting
$\rho_{i,\ell_k}:=|\alpha_{\ell_k+1}|\,|(U_B)_{\ell_k+1,i}|$ the residual bound of the $i$-th Ritz triplet, since $\tilde\sigma_{i,\ell_k}(M_{\ell_k})\le\sigma_i(\widetilde Y^{k-1})$, we have $0\le\sigma_i(\widetilde Y^{k-1})-\tilde\sigma_{i,\ell_k}(M_{\ell_k})\le\rho_{i,\ell_k}$ by (\ref{eq:Lanczos_error_bound_v2}) and
{\small
\begin{align*}
\sigma_i(\widetilde Y^{k-1})^2-\tilde\sigma_{i,\ell_k}(M_{\ell_k})^2=&\,(\sigma_i(\widetilde Y^{k-1})-\tilde\sigma_{i,\ell_k}(M_{\ell_k}))(2\tilde\sigma_{i,\ell_k}(M_{\ell_k})+\sigma_i(\widetilde Y^{k-1})-\tilde\sigma_{i,\ell_k}(M_{\ell_k}))\\
\le&\,\rho_{i,\ell_k}(2\tilde\sigma_{i,\ell_k}(M_{\ell_k})+\rho_{i,\ell_k}),
\end{align*}
}
because $t\mapsto t(2\tilde\sigma_{i,\ell_k}+t)$ is nondecreasing for $t\ge0$. By~\eqref{eq:error_lanczos},
we get
\begin{align*}
\|\Delta^k\|_F^2\le&\;\|M_{\ell_k}-\widetilde Y^{k-1}\|_F^2-\sum_{i>\ell_k}\min(\sigma_i(\widetilde Y^{k-1}),\tau)^2\\
=&\;\|M_{\ell_k}\|_F^2+\|\widetilde Y^{k-1}\|_F^2-2\text{trace}(M_{\ell_k}^T\widetilde Y^{k-1})-\sum_{i>\ell_k}\sigma_i(\widetilde Y^{k-1})^2,
\end{align*}
where we exploit the assumption that $\sigma_i(\widetilde Y^{k-1})\leq \tau$ for any $i>\ell_k$. Now, we observe that by exploiting the cyclic property of the trace operator and the Lanczos recurrence~\eqref{eq:rel_lanczos2}, we can write
\begin{align*}
\text{trace}(M_{\ell_k}^T\widetilde Y^{k-1})=&\;\text{trace}(Q_{\ell_k}B_{\ell_k}^TP_{\ell_k+1}^T\widetilde Y^{k-1})=\text{trace}(B_{\ell_k}^TP_{\ell_k+1}^T\widetilde Y^{k-1}Q_{\ell_k})\\
=&\;\text{trace}(B_{\ell_k}^TB_{\ell_k})=\sum_{i=1}^{\ell_k}\tilde\sigma_{i,\ell_k}^2=\|M_{\ell_k}\|_F^2.
\end{align*}
Therefore, 
\begin{equation}\label{eq:bnd_full}
\|\Delta^k\|_F^2\leq 
\sum_{i\le \ell_k}\Big(\sigma_i(\widetilde Y^{k-1})^2-\tilde\sigma_{i,\ell_k}(M_{\ell_k})^2\Big) \le \sum_{i=1}^{\ell_k}\rho_{i,\ell_k}\,(2\tilde\sigma_{i,\ell_k}(M_{\ell_k})+\rho_{i,\ell_k})=:\tilde\omega_{\ell_k}.
\end{equation}
Stopping when $\tilde\omega_{\ell_k}<C/k^{2+\epsilon}$ guarantees 
the bounds on $\|\Delta^k\|_F$ and $\|\Delta^k\|_*$,  so that the convergence
theory applies. Note that $\tilde\omega_{\ell_k}$ involves \emph{all} the $\ell_k$ Ritz values, including
those below $\tau$, even though the latter ones will be then discarded by the soft-thresholding operator. The next heuristic criterion builds upon this observation.

\paragraph{Heuristic  bound - \frobbndrk}

Let $r_k$ denote the number of Ritz values larger than the threshold $\tau$, i.e.
$r_k=\mathrm{rank}\big(\mathcal S_\tau(M_{\ell_k})\big)$. Then we can restrict
the sum in~\eqref{eq:bnd_full} to the $r_k$ Ritz values above $\tau$, define
the quantity 
\begin{equation}\label{eq:bnd_rk}
\bar \omega_{\ell_k} := \sum_{i\le r_k}\rho_{i,\ell_k}(2\tilde\sigma_{i,\ell_k}(M_{\ell_k})+\rho_{i,\ell_k}),
\end{equation}
and stop when $\bar \omega_{\ell_k}<\frac{C}{k^{2+\epsilon}}$.
To the best of our knowledge, this cannot be seen as a rigorous bound for $\|\Delta^k\|_F^2$: the terms
with $\tilde\sigma_{i,\ell_k}(M_{\ell_k})\le\tau$ are neglected, although in general they do not vanish.
\vspace{0.2cm}

Table \ref{tab:criteri_riassunto} summarized the inexact stopping criteria for the Lanczos bidiagonalization to approximate the SVD of $\widetilde Y^{k-1}$. It also includes the cost of evaluating the stopping tests and their relations with the convergence theory.
\begin{table}[t]
\centering
\small
\renewcommand{\arraystretch}{1.4}
\begin{tabular}{@{}l r>{\raggedright\arraybackslash}p{3cm} r r@{}}
\toprule
Criterion & Bound  & Guarantee on $\|\Delta^k\|_F^2$ & Cost (flops)\\
\midrule
\frob
  & $\omega_{\ell_k}=\displaystyle \|M_{\ell_k}-\widetilde Y^{k-1}\|_F^2$
  & rigorous bound
  & $\mathcal O(1)$ \\[6pt]
\frobbnd
  & $\tilde \omega_{\ell_k}=\displaystyle \sum_{i=1}^{\ell_k}\rho_{i,\ell_k}\,(2\tilde\sigma_{i,\ell_k}(M_{\ell_k})+\rho_{i,\ell_k})$
  & rigorous bound (w/ some assumptions)
  & $\mathcal O(\ell_k)$ \\[6pt]
\frobbndrk
  & $\bar \omega_{\ell_k}=\displaystyle \sum_{i=1}^{r_k}\rho_{i,\ell_k}\,(2\tilde\sigma_{i,\ell_k}(M_{\ell_k})+\rho_{i,\ell_k})$
  & heuristic
  & $\mathcal O(\ell_k)$ \\[6pt]

\bottomrule
\end{tabular}
\caption{Stopping criteria for the Lanczos bidiagonalization. Here $\rho_{i,\ell_k}=\alpha_{\ell_k+1}|(U_B)_{\ell_k+1,i}|$ is the residual bound of the $i$-th Ritz value $\tilde\sigma_{i,\ell_k}(M_{\ell_k})$, and $r_k=|\{i\le\ell_k:\tilde\sigma_{i,\ell_k}(M_{\ell_k})>\tau\}|$.
All the criteria stop Lanczos at the first step at which the quantity in the second column is smaller than $C/k^{2+\epsilon}$.}
\label{tab:criteri_riassunto}
\end{table}

\section{Numerical experiments}\label{sec_Numerical_results}
In this section we report some numerical results to support the theoretical findings derived in the previous sections. To this end, we considered the Matlab implementation of SVT as described in Algorithm \ref{alg:svt}, 
available at \url{https://candes.su.domains/software/svt/code.html} and conducted all experiments in MATLAB R2026a on an Intel Ultra 9 processor with 128GB DDR5 RAM. 
In all our numerical experiments, the sequence of stepsizes $\{\delta_k\}$ is chosen constant (parameter $\delta$ in the following) and the starting point is chosen as  $\widetilde Y^0=k_0\delta P_\Omega(A)$, 
$k_0=\lceil\tau/(\delta\|P_\Omega(M)\|_2)\rceil$; see \cite{CandesRecht2009}. 
The matrix $\widetilde Y^{k-1}$  is stored as a sparse matrix and accessed only through matrix-vector products. At each SVT iteration we employ the Lanczos bidiagonalization with partial reorthogonalization of PROPACK (\texttt{lansvd}, \texttt{lanbpro}) for approximating the SVD of $\widetilde Y^{k-1}$; see \cite{propack_doc,larsen2004propack}. These routines are used as a black box except for the stopping test that we changed according to our analysis. If $r_{k-1}:=\text{rank}(\mathcal{S}_\tau(\widetilde Y^{k-2}))$, at the $k$th SVT iteration we ask 
PROPACK to compute up to $s_k=r_{k-1}+1$ singular triplets of $\widetilde Y^{k-1}$.
If the $s_k$th computed singular value is still larger than $\tau$, some relevant singular values might have been missed.  In this case, the number of requested triplets is increased by a parameter $b$ and the computation is repeated,  until at least one singular value smaller than or equal to $\tau$ is found, thus reducing the chances of discarding significant components while keeping the computational cost minimal.

In the original SVT implementation, the default setting makes PROPACK stop the Lanczos process as soon as the residual norm related to the requested $s_k$ Ritz triplets is sufficiently small in a relative sense, that is, when
\begin{equation}\label{eq:svt_test}
    \max_{i=1,\dots,s_k}\tilde \rho_{i,s_k}/\tilde\sigma_{i,s_k}(M_{\ell_k})\le 16\,\varepsilon_{\rm mach},
\end{equation}
where $\tilde \rho_{i,s_k}$ denotes the upper bound given by the PROPACK routine \texttt{refine\-bounds} \cite[Appendix B]{propack_doc} on the residual norm achieved by the $i$-th approximate singular triplet. We check this condition at every step $\ell_k$, since its cost is negligible with respect to that of matrix--vector products with $\widetilde Y^{k-1}$\footnote{In the default setting, PROPACK checks~\eqref{eq:svt_test} only at the end of a number of Lanczos steps.}. Although extremely effective in practice as illustrated in the following, to the best of our knowledge the condition~\eqref{eq:svt_test} does not guarantee a decrease in $\|\Delta^k\|_F^2$ and  $\|\Delta^k\|_*^2$ as required by our convergence analysis of SVT. In this regard, also~\eqref{eq:svt_test} must be viewed as a heuristic stopping criterion. In the results that follow we will denote the criterion~\eqref{eq:svt_test} by \orig. 

Concerning the stopping criteria proposed in this paper, we monitor at every Lanczos step $\ell_k$ one of the following quantities: $\omega_{\ell_k}$ in~\eqref{eq:Lancsoz_stop} (\frob), $\tilde\omega_{\ell_k}$ in~\eqref{eq:bnd_full} (\frobbnd), or
$\bar\omega_{\ell_k}$ in~\eqref{eq:bnd_rk} (\frobbndrk). Lanczos is stopped as soon as this quantity becomes smaller than $C/k^{2+\epsilon}$ for given $C$ and $\epsilon$. Following what we do with these criteria, we will also consider a relaxed variant of~\eqref{eq:svt_test}. In particular, the condition
\begin{equation}\label{eq:svt_test_abs}
    \max_{i=1,\dots,s_k}\tilde\rho_{i,s_k}^{\,2}\;\le\;\frac{C}{k^{2+\epsilon}} ,
\end{equation}
will be monitored; we will refer to this stopping criterion as \origtol. Notice that the left-hand side in~\eqref{eq:svt_test_abs} is now an absolute quantity in contrast to what happens in~\eqref{eq:svt_test}. This is to make any comparison between \origtol\ and the absolute conditions in \frob, \frobbnd, and \frobbndrk as fair as possible. We mention once again that the quantity monitored by \origtol\ (and also \orig) is not a bound on the error $\|\Delta^k\|_F^2$ , since it only controls the residuals of the single Ritz triplets, so that employing this criterion does not guarantee the convergence of the overall SVT iteration.

All the variants of the Lanczos scheme use the same starting vector (fixed random seed). Moreover, in case the full ambient space is generated, Lanczos gets stopped even though the considered stopping criterion may not be satisfied.

\paragraph{Parameters}
We use $b=5$, $\epsilon=10^{-6}$, and $C$ was varied in the experiments.
A very small $\epsilon$ makes the thresholds $C/k^{2+\epsilon}$ summable in $k$ (so the errors of the successive iterations cannot accumulate) while decaying as slowly as possible, which allows Lanczos to stop early. The value $\delta=1.2\,n_1n_2/|\Omega|$ is the standard one in the original SVT algorithm.
The maximum iteration count is set equal to $k_{max}= 1000$.

\subsection{Preliminary tests}
In this section we compare the original SVT (\orig) with the inexact variants as described in the previous section.

\begin{figure}[t]
\centering
\includegraphics[width=\linewidth]{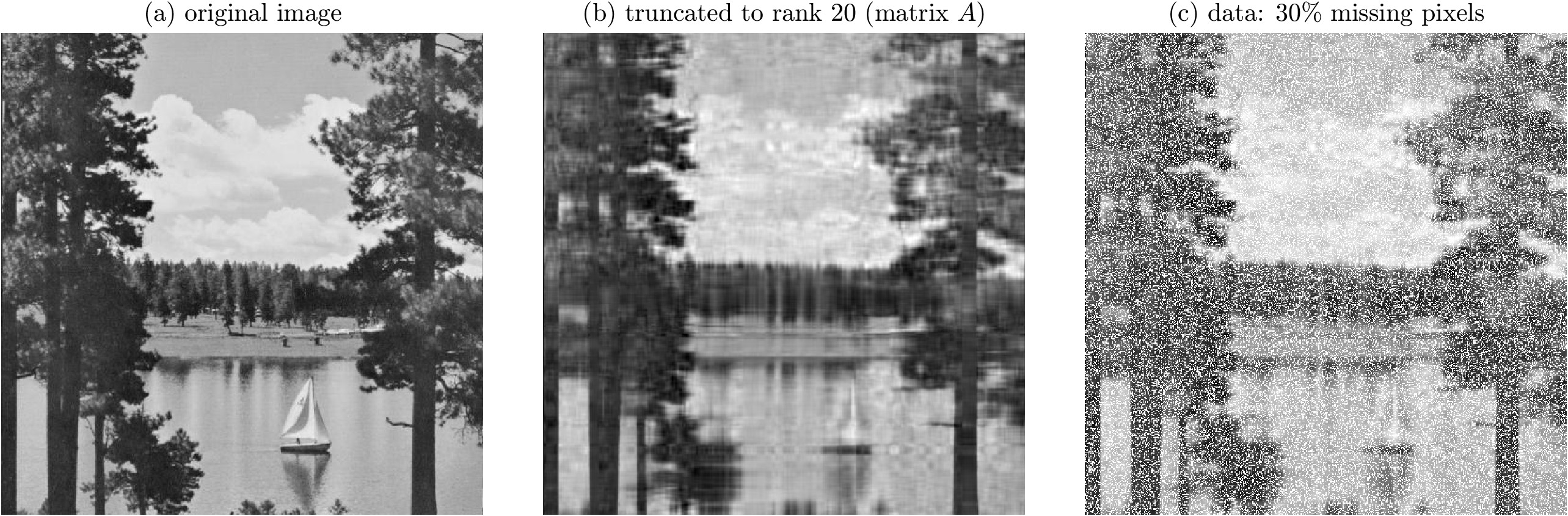}
\caption{Test problem \texttt{lake}: (a) original image $512\times512$; (b) the same image truncated to rank $20$ by SVD; this represents the matrix $A$ to be reconstructed; (c) data given to SVT, with $30\%$ of randomly chosen missing pixels (in white).}
\label{fig:lake_problema}
\end{figure}

As test problem, we consider the completion of a $512\times512$ grayscale image (\texttt{lake}), stored as a matrix of intensities.
The matrix is first truncated to rank $r=20$ by its singular value decomposition, and the truncated matrix is the matrix $A$ to be reconstructed.
A set $\Omega$ of observed entries is obtained by removing $30\%$ of the pixels, chosen uniformly at random, so that $70\%$ of the entries are observed; see Figure \ref{fig:lake_problema}.

We used the threshold $\tau=10^4$ leading to a solution of rank $20$ in all tests and $tol = 10^{-3}$ in the  SVT algorithm. 
 The latter equipped with \orig needs $24$ iterations to converge, $1461$ total Lanczos steps, and reaches a relative error of $1.18\cdot10^{-3}$. We will consider these results as baseline to compare the performance of the inexact criteria with. In Table~\ref{tab:lake_vary_C} we report the results achieved by 
 \frob, \frobbnd, \frobbndrk, and \origtol\ by using $C=10^2,10^5$.
All the inexact criteria, except for \origtol, need the same $24$ SVT iterations to converge and reach the same final error, but they differ considerably in the number of Lanczos steps.
In particular, \frobbndrk\ manages to save $16\%$ ($C=10^2)$ to almost $22\%$ ($C=10^5$) Lanczos steps compared to \orig, whereas \frob\ and \frobbnd\ require significantly  more steps.
This is due to the fact that these criteria are too conservative and the quantities $\omega_{\ell_k}$ and $\tilde\omega_{\ell_k}$ do not decrease below the given threshold until the Lanczos process has generated (almost) the whole space, whose dimension is $512$, leading to a tremendous number of Lanczos steps overall.
It is interesting to observe how poorly \origtol performs. Similarly to \frobbndrk, this is an heuristic criterion. However, Table~\ref{tab:lake_vary_C} shows that the latter is much more reliable than the former. We further investigate this aspect in the following experiment.

\begin{table}[t]
\centering
\small
\begin{tabular}{lccc}
\multicolumn{4}{l}{\emph{\texttt{lake} $512\times512$, $30\%$ missing pixels}} \\
\toprule
 & SVT it. & $\ell_{\rm tot}$ & ${\|X-A\|_F}/{\|A\|_F}$ \\
\midrule

\orig & 24 & 1461 & $1.18\cdot10^{-3}$ \\
\cmidrule(l){1-4}
\multicolumn{4}{l}{\emph{\origtol}} \\
$C=10^{2}$ & 345 & 14792\ (-912.5\%) & $8.31\cdot10^{-4}$ \\
$C=10^{5}$ & 1000 & 21695\ (-1384.9\%) & $1.06\cdot10^{-1}$ \\
\cmidrule(l){1-4}
\multicolumn{4}{l}{\emph{\frob}} \\
$C=10^{2}$ & 24 & 15360\ (-951.3\%) & $1.18\cdot10^{-3}$ \\
$C=10^{5}$ & 24 & 15337\ (-949.8\%) & $1.18\cdot10^{-3}$ \\
\cmidrule(l){1-4}
\multicolumn{4}{l}{\emph{\frobbnd}} \\
$C=10^{2}$ & 24 & 15360\ (-951.3\%) & $1.18\cdot10^{-3}$ \\
$C=10^{5}$ & 24 & 15354\ (-950.9\%) & $1.18\cdot10^{-3}$ \\
\cmidrule(l){1-4}
\multicolumn{4}{l}{\emph{\frobbndrk}} \\
$C=10^{2}$ & 24 & 1227\ (16.0\%) & $1.18\cdot10^{-3}$ \\
$C=10^{5}$ & 24 & 1142\ (21.8\%) & $1.18\cdot10^{-3}$ \\

\bottomrule
\end{tabular}
\caption{Test problem: image \texttt{lake} $512\times512$ truncated to rank $20$ with $30\%$ of missing pixels. SVT with the original test (\orig) and SVT with \frob, \frobbnd, \frobbndrk\ and with the original test with inexact tolerance (\origtol) for increasing values of $C$ in the threshold $C/k^{2+\epsilon}$,  $\tau=10^4$. $\ell_{\rm tot}$ is the total number of Lanczos steps; in parentheses the saving with respect to \orig\ (negative if more steps are needed).}
\label{tab:lake_vary_C}
\end{table}

\begin{figure}[t]
\centering
\includegraphics[width=\linewidth]{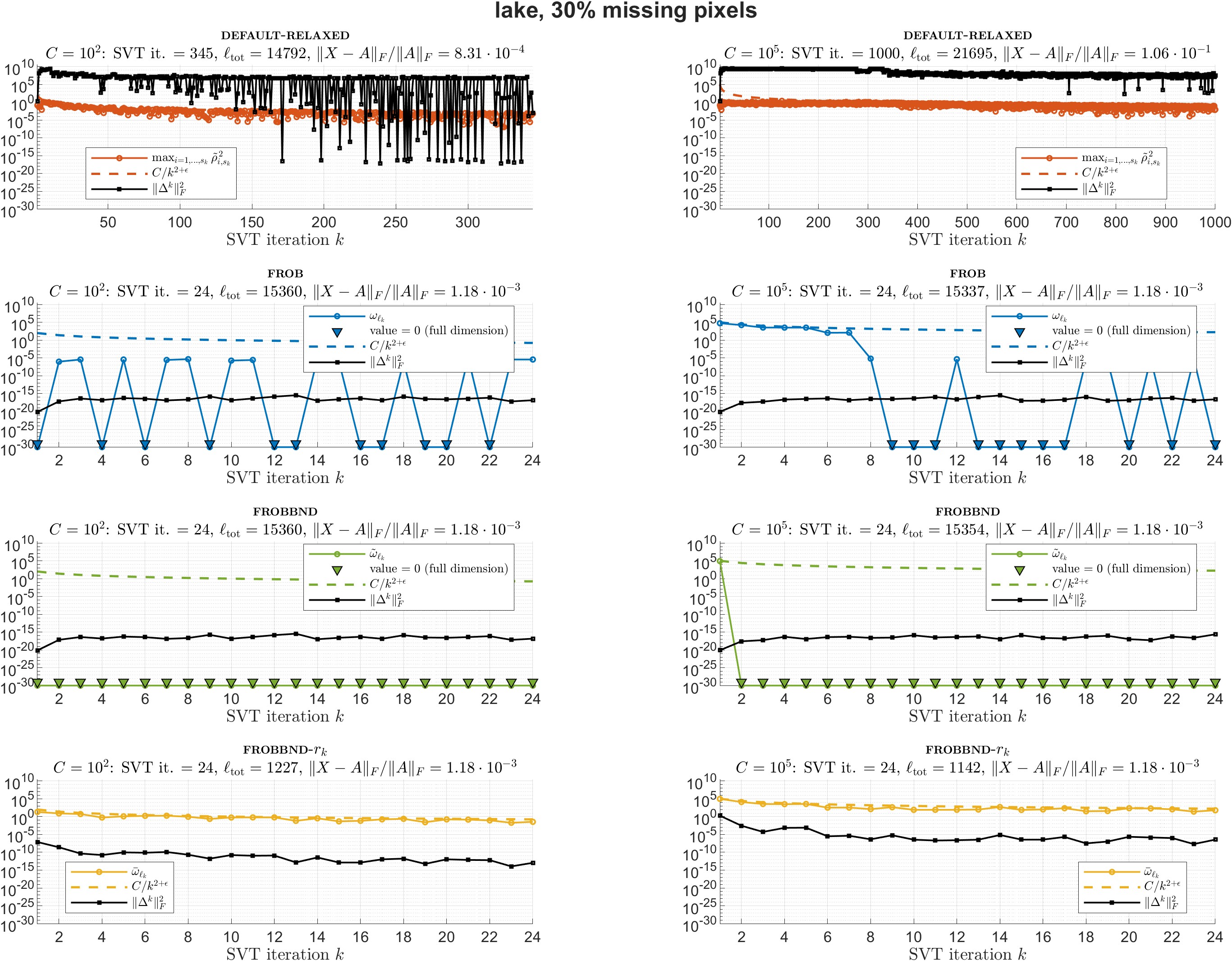}
\caption{SVT on the image \texttt{lake} $512\times512$ truncated to rank $20$ with $30\%$ of missing pixels, with the four criteria: from top to bottom \origtol, \frob, \frobbnd, and \frobbndrk; left column $C=10^{2}$, right column $C=10^{5}$. For each SVT iteration $k$, evaluated at the end of the last Lanczos call of the iteration, each panel reports the quantity monitored by the criterion (solid line), its threshold $C/k^{2+\epsilon}$ (dashed line), and the true error $\|\Delta^k\|_F^2$ of the same run. 
The triangles at the bottom of the panels mark the iterations in which the monitored quantity is zero (the whole space has been generated).}
\label{fig:lake_r20_4x2_C100_C100000}

\end{figure}

Figure~\ref{fig:lake_r20_4x2_C100_C100000} shows the results obtained by using  \origtol\ (first from top), \frob\ (second from top), \frobbnd\ (third from top), and \frobbndrk\ (bottom) for $C=10^2$ (left) and $C=10^5$ (right).
Each panel reports the number of SVT iterations, the total number $\ell_{\rm tot}$ of Lanczos steps performed by the overall SVT run, and the relative error $\|X-A\|_F/\|A\|_F$ of the final iterate.
The curves are functions of the SVT iteration $k$ and refer to the end of the last Lanczos call of each iteration: the quantity monitored by the criterion (solid line), its threshold (dashed line), and the true error $\|\Delta^k\|_F^2$ in (\ref{eq:error_lanczos}), where $\mathcal S_\tau(\widetilde Y^{k-1})$ is obtained from the exact singular value decomposition of $\widetilde Y^{k-1}$.
The results in Figure~\ref{fig:lake_r20_4x2_C100_C100000} (very top) explain the failure of~\origtol. Indeed, in this case, $\|\Delta_k\|_F^2$ (black solid line) turns out to be always larger than $C/k^{2+\epsilon}$ so that the results derived in section~\ref{sec_inexactSVT} cannot be applied to ensure the convergence of SVT and the latter does indeed perform rather poorly.
Similarly, the second and third panels from the top of Figure~\ref{fig:lake_r20_4x2_C100_C100000} illustrate how conservative \frob\ and \frobbnd are. Indeed, the triangles in these panels indicate at which SVT iteration Lanczos had to build the whole subspace to meet the stopping criteria given by \frob\ and \frobbnd. For all these iterations, $\|M_{\ell_k}-\widetilde{Y}^{k-1}\|_F^2=0$ so that $\|\Delta^{k}\|_F^2$ matches machine precision. The condition \frobbndrk\ is the only one that really acts like an ``inexact'' criterion. It allows to save a sizable number of Lanczos iterations at the beginning of the SVT process to then becoming stricter for larger $k$. In spite of this inexactness, $\|\Delta_k\|_F^2$ remains well below $C/k^{2+\epsilon}$ as shown in Figure~\ref{fig:lake_r20_4x2_C100_C100000} (bottom) and the overall SVT scheme thus converges.

To recap, these experiments show that the criteria \frob\ and \frobbnd\ are too conservative: they stop the Lanczos process only when (almost) the whole space has been generated, so they give no saving with respect to \orig.
The criterion \origtol\ does not provide a bound on the error (neither in theory nor in practice), and therefore it has no guarantee of convergence.
Although still heuristic, \frobbndrk provides a practical bound for $\|\Delta_k\|_F^2$ and allows an implementation in which an inexact approach to SVT can be tested. This last point is discussed in more detail in the following section.

\subsection{Further experiments with \frobbndrk}

In this section, we challenge the robustness of the heuristic \frobbndrk as this turned out to be the most performing criterion in the previous section. We therefore consider a set of random low-rank matrix completion problems where $A\in\mathbb R^{n\times m}$ of rank $r$ is generated as $A=GH$, where $G\in\mathbb R^{n\times r}$ and $H\in\mathbb R^{r\times m}$ have independent standard Gaussian entries. 

We report two sets of experiments both comparing SVT with \orig and \frobbndrk\ where for the latter we vary the value of $C\in\{10,\ldots,10^6\}$.
In the first set of experiments (Table~\ref{tab:random_grandi_frobbndrk_C_ranks}) we study the dependence on the rank and on the size: we assume $60\%$ of the entries of $A$ are known and consider the ranks $r=10,30$, and the sizes  $2000\times600$, $5000\times1000$, and $10000\times2000$.
In the second one (Table~\ref{tab:random_grandi_frobbndrk_C_ranks_obs30_dfac12}) we reduce the number of known observations to $30\%$ but keep the same rank values and dimensions. 
For each run we report the CPU time in seconds, the number of SVT iterations, the total number $\ell_{\rm tot}$ of Lanczos steps performed in the whole SVT run, and the relative error $\|X-A\|_F/\|A\|_F$ of the final iterate. The adopted value of  $\tau$ is specified in the table of results and is chosen empirically in order to recover a rank $r$ matrix.

\input{tabelle_det/tabella_random_grandi_frobbndrk_C_ranghi10e30_obs60}
\input{tabelle_det/tabella_random_grandi_frobbndrk_C_ranghi10e30_obs30_dfac12}

Some comments are as follows. For $C\le10^{5}$, \frobbndrk consistently leads to the execution of the same number of iterations (up to a single unit) SVT performs when equipped with \orig. Similarly, the accuracy of the final iterate is comparable to that achieved when employing \orig.  For $C\le10^{5}$, \frobbndrk is able to lower the overall number of Lanczos steps and, consequently, the running times. These savings are substantial in most of the tests, up to about one half of the Lanczos steps and of the CPU time. It is interesting to notice that the major gains come for $A$ with a smaller rank and larger dimensions, a sizable percentage of known entries, and large $C$. We must mention, though, that sensible selections of the latter parameter must be adopted. Indeed, in all our tests, selecting $C=10^{6}$ significantly deteriorated the results achieved by SVT with \frobbndrk. Indeed, this algorithm needs more SVT iterations than \orig in general, leading, in some cases, to performing even more Lanczos steps overall.

    \section{Conclusions}\label{sec_Conclusions}

Adopting iterative procedures for SVD approximation is often the only viable option to achieve competitive performance of the SVT algorithm, especially for large-scale problems. This common practice introduces some inexactness in the latter iterative procedure that, to our knowledge, has never been accounted for in the study of the convergence of SVT. This paper fills this gap by deriving the right theoretical framework able to guarantee the convergence of SVT also in case of inexact SVD computations. Our results have been then specialized to the case of Lanczos bidiagonalization for SVD approximation. In particular, we showed how to introduce cheap stopping criteria in the Lanczos scheme able to guarantee the overall convergence of the underlying SVT process. Our convergence analysis also suggested the introduction of sensible heuristic criteria that turned out to be reliable and competitive in practice.

While the analysis performed in this paper envisions the adoption of a broad class of well-established methods for SVD approximation, it is not able to include emerging randomized schemes like randomized SVD~\cite{halko2011finding}. The convergence analysis of SVT equipped with randomization-based algorithms for SVD approximation requires to resort to stochastic error analysis; 
a completely different approach than the one followed here. We plan to explore this challenging venue in the near future.


{\footnotesize

\subsection*{\footnotesize Ethics approval and consent to participate}
Not applicable.
\subsection*{\footnotesize Consent for publication}
Not applicable.
\subsection*{\footnotesize Funding}
{\footnotesize
The work of L.B., D.P., and M.P. was partially supported by INdAM-GNCS under the
project CUP  E53C25002010001.  The research of M.P.  was  partially granted by the Italian Ministry of University and Research (MUR) through the PRIN 2022 ``MOLE: Manifold constrained Optimization and LEarning'',  code: 2022ZK5ME7 MUR D.D. financing decree n. 20428 of November 6th, 2024 (CUP B53C24006410006). 

}
\subsection*{\footnotesize Availability of data and materials}
The data presented in this work are available from the corresponding author upon request.

\subsection*{\footnotesize Competing interests}
The authors declare that they have no competing interests.
\subsection*{\footnotesize Authors' contributions}
All authors contributed equally to the writing of this article. All authors reviewed the manuscript.}

\subsection*{\footnotesize AI declaration}
The research direction, the method and the approach to the theoretical analysis, including proofs, 
were developed by the authors without the employment of any AI tool. Claude Sonnet 5.5 was prompted 
for the copy-editing the text and to implement the inexact stopping criteria within the original SVT code available at \url{https://candes.su.domains/software/svt/code.html},
prepare the benchmarks, executables for hyperparameter studies and with preparing the plots for the paper.

\section*{Appendix} \label{app:svt}

To make this work as self-contained as possible, in this section we report some technical results, with the accompanying proof, useful for our analysis.

\begin{lemma}\label{stima_subdiff}
If $Z\in \partial\|X\|_*$ then $\|Z\|_2\leq 2$.
\end{lemma}
\begin{proof}
    As stated in the article \cite{halko2011finding}, the subdifferential of the nuclear norm at $X=U\Sigma V^T$ can be characterized as follows
    \begin{equation}\label{nuclear_subgradient}
        \partial \|X\|_*=\{UV^T+W|W\in \mathbb{R}^{n_1\times n_2}, U^TW=0, WV=0, \|W\|_2\leq 1\}.
    \end{equation}
For every $Z\in \partial \|X\|_*$, some $W$ exists such that $Z=U^TV+W$, hence
\begin{equation*}
    \|Z\|_2\leq\|U^TV\|_2+\|W\|_2\leq \|U^T\|_2\|V\|_2+\|W\|_2\leq 2. 
\end{equation*}
For the last inequality we use that the matrices $U$ and $V$ are orthonormal and the characterization of the matrix $W$ given in \eqref{nuclear_subgradient}.
\end{proof}

{
\begin{lemma}\label{lemma_char_Z}
Let $f_\tau$ be defined in \eqref{PB_SVT}, then 
 $Z \in \partial f_\tau(X)$ if and only if $Z = \tau Z_0 + X$, where $Z_0 \in \partial X^*$.
\end{lemma}
\begin{proof}
    The claim follows from the sum rule for subdifferentials, along with the fact that $\tfrac12\|X\|_F^2$ is differentiable with gradient $X$.
    \end{proof}
}

\begin{definition}[Two-norm ball Euclidean projection]
Let $Y \in \mathbb{R}^{n_1 \times n_2}$ be a matrix with singular value decomposition $Y = U \Sigma V^T$,
let $p = \min\{n_1, n_2\}$ and let $\sigma_1(Y) \ge \dots \ge \sigma_p (Y)\ge 0$ be the singular values of $Y$. Let a parameter $\tau > 0$ be fixed. The \textit{two-norm ball Euclidean projection} $\mathcal{P}_\tau$ is an operator from $\mathbb{R}^{n_1 \times n_2}$ to $\mathbb{R}^{n_1 \times n_2}$ defined as
\begin{equation}\label{Operator_shrinkage_bis}
    \mathcal{P}_\tau(Y) = U \mathcal{P}_\tau(\Sigma) V^T, \qquad
    (\mathcal{P}_\tau(\Sigma) )_{ii}=  \min(\sigma_i(Y),\tau), \quad i=1,\dots,p.
\end{equation}
\end{definition}
\begin{remark}
For all  matrix $Y$ we have
    \begin{equation*}
        Y-\mathcal S_\tau(Y)=\mathcal{P}_\tau(Y).
    \end{equation*}
Indeed, let $r$ be the last index such that $\sigma_r(Y)>\tau$ and $p = \operatorname{min}\{n_1, n_2\}$. It holds $Y-\mathcal S_\tau(Y)=U\big(\Sigma-\mathcal S_\tau(\Sigma)\big)V^T$ and $\Sigma-S_\tau(\Sigma)=\operatorname{diag}\big(\min(\sigma_1(Y),\tau), \dots,$ $\ min(\sigma_p(Y),\tau)\big)$.
    Hence, 
    \begin{equation*}
        Y-\mathcal S_\tau(Y)=U\mathcal{P}_\tau(\Sigma)V^T=\mathcal{P}_\tau(Y).
    \end{equation*}
\end{remark}

In the following we provide and revisit some proofs of the paper~\cite{oh2018fast}. In particular the following Lemma is an extended version of the Proposition 2 of~\cite{oh2018fast}.
\begin{lemma}\label{est_op_shr_fr}
Let $\mathcal S_\tau(\cdot)$ be the shrinkage operator; then the following inequality holds for any $ A,B\in \mathbb{R}^{n_1\times n_2}$
\begin{equation}\label{pseudo-contraction}
    \|\mathcal S_\tau(A)-\mathcal S_\tau(B)\|_F^2\leq \|A-B\|_F^2-\|\mathcal P_\tau(A)-\mathcal P_\tau(B)\|_F^2.
\end{equation}
\end{lemma}
\begin{proof}
Let us consider
    \begin{align*}
        \|A-B\|_F^2-\|\mathcal{P}_\tau(A)-\mathcal{P}_\tau(B)\|_F^2=&\|\mathcal S_\tau(A)+\mathcal{P}_\tau(A)-(\mathcal S_\tau(B)+\mathcal{P}_\tau(B))\|_F^2\\
        &-\|\mathcal{P}_\tau(A)-\mathcal{P}_\tau(B)\|_F^2\\
        =&\|\mathcal S_\tau(A)-\mathcal S_\tau(B)\|_F^2 \\
        & + 2\langle \mathcal S_\tau(A)-\mathcal S_\tau(B),\mathcal{P}_\tau(A)-\mathcal{P}_\tau(B)\rangle.
    \end{align*}
    Hence, to proof the inequality \eqref{pseudo-contraction} we prove $0\leq\langle \mathcal S_\tau(A)-\mathcal S_\tau(B),\mathcal{P}_\tau(A)-\mathcal{P}_\tau(B)\rangle$. We have
    \begin{align*}
        \langle \mathcal S_\tau(A)-\mathcal S_\tau(B),\mathcal{P}_\tau(A)-\mathcal{P}_\tau(B)\rangle&=\langle \mathcal S_\tau(A)-\mathcal S_\tau(B),A-S_\tau(A)-B+S_\tau(B)\rangle\\
        =&-\|\mathcal S_\tau(A)-\mathcal S_\tau(B)\|_F^2 \\
        &+\langle \mathcal S_\tau(A)-\mathcal S_\tau(B),A-B\rangle.
    \end{align*}
    We have to prove that $\|\mathcal S_\tau(A)-\mathcal S_\tau(B)\|_F^2\leq\langle \mathcal S_\tau(A)-\mathcal S_\tau(B),A-B\rangle$.
Recall that $\mathcal S_\tau(A)$ is the unique minimizer of
$f^A_\tau(X):=\tau\|X\|_*+\frac12\|X-A\|_F^2$ \cite{cai2010singular}, calculating the gradient we obtain
\begin{align*}
    0&= \mathcal S_\tau(A)-A+\tau Z_A, \ \text{with } Z_A\in \partial \|\mathcal S_\tau(A)\|_*,\\
    0&= \mathcal S_\tau(B)-B+\tau Z_B, \ \text{with } Z_B\in \partial \|\mathcal S_\tau(B)\|_*.
\end{align*}
Hence,
\begin{equation*}
    \mathcal S_\tau(A)-A-(B-\mathcal S_\tau(B))=\tau(Z_A-Z_B).
\end{equation*}
Considering the scalar product with $\mathcal S_\tau(A)-\mathcal S_\tau(B)$ we have
\begin{equation*}
    \langle A-B,\mathcal S_\tau(A)-\mathcal S_\tau(B)\rangle-\|\mathcal S_\tau(A)-\mathcal S_\tau(B)\|_F^2=\tau\langle Z_A-Z_B,\mathcal S_\tau(A)-\mathcal S_\tau(B)\rangle\geq0.
\end{equation*}
The last inequality is true for \cite[Lemma 4.1]{cai2010singular}.
This proves the lemma.
\end{proof}

The proof of the following lemma is a revised version of the proof of Lemma 2 in \cite{oh2018fast}.
\begin{lemma}\label{error_bound_Ptau}
Let $Y, A\in \mathbb{R}^{n_1\times n_2}$ be matrices such that $A$ has rank-$p$. Then
\begin{equation}\label{stima_op_p}
    \min \|\mathcal{P}_\tau(Y)-\mathcal{P}_\tau(A)\|_F^2=\sum_{i>p}\min(\sigma_i(Y),\tau)^2
\end{equation}
\end{lemma}
\begin{proof}
Notice that if $A$ has rank $p$, also $\mathcal P_\tau(A)$ has rank $p$. Then, using the unitary invariant of the norm, 
\begin{align*}
    \|\mathcal P_\tau(Y)-\mathcal P_\tau(A)\|^2_F
    &=\|U_Y\mathcal P_\tau(\Sigma_Y)V_Y^T-\mathcal P_\tau(A)\|^2_F= \|\mathcal P_\tau(\Sigma_Y)-U_Y^T\mathcal P_\tau(A)V_Y\|^2_F.
\end{align*}
Now we can notice that the matrix $U_Y^T\mathcal P_\tau(A)V_Y$ has rank $p$, hence the minimum error in Frobenius norm is given considering the truncated rank-$p$ SVD of $\mathcal P_\tau(\Sigma_Y)$, which is
$\mathcal P_\tau(\Sigma_Y)_p= \operatorname{diag}\big(\min(\sigma_1(Y),\tau),\dots,\min(\sigma_p(Y),\tau),0,\dots,0\big).$
Hence, we have, defining $n=\min(n_1,n_2)$
\begin{align*}
\mathcal     P_\tau(\Sigma_Y)-\mathcal P_\tau(\Sigma_Y)_p&=\operatorname{diag}\big(0,\dots,0,\min(\sigma_{p+1}(Y),\tau),\dots,\min(\sigma_n(Y),\tau)\big),
\end{align*}
so that
    \begin{align*}
        \|\mathcal P_\tau(Y)-\mathcal P_\tau(A)\|^2_F&= \|\mathcal P_\tau(\Sigma_Y)-U_Y^T\mathcal P_\tau(A)V_Y\|^2_F\geq \sum_{i>p}\min(\sigma_i(Y),\tau)^2.
    \end{align*}
\end{proof}

\bibliographystyle{spmpsci}

\end{document}

%% file: tabelle_det/tabella_random_grandi_frobbndrk_C_ranghi10e30_obs60.tex
\begin{table}[t]
\centering
\small
\begin{tabular}{lcccc}
\toprule
 & Time (s) & SVT it. & $\ell_{\rm tot}$ & 
 $\|X-A\|_F/\|A\|_F$ \\
\midrule
\multicolumn{5}{l}{\emph{Random $2000\times600$, rank $10$, $60\%$ observed, $\tau=1500$}} \\
\orig & 1.36 & 27 & 935 & 
$9.08\cdot10^{-7}$ \\
\cmidrule(l){1-5}
\frobbndrk, $C=10$ & 0.82 & 27 & 641 & 
$9.08\cdot10^{-7}$ \\
\frobbndrk, $C=10^{2}$ & 0.81 & 27 & 629 & 
$9.08\cdot10^{-7}$ \\
\frobbndrk, $C=10^{3}$ & 0.78 & 27 & 612 & 
$9.08\cdot10^{-7}$ \\
\frobbndrk, $C=10^{4}$ & 0.76 & 27 & 597 & 
$9.15\cdot10^{-7}$ \\
\frobbndrk, $C=10^{5}$ & 0.76 & 27 & 592 & 
$8.94\cdot10^{-7}$ \\
\frobbndrk, $C=10^{6}$ & 1.37 & 52 & 1079 & 
$9.71\cdot10^{-7}$ \\
\midrule
\multicolumn{5}{l}{\emph{Random $2000\times600$, rank $30$, $60\%$ observed, $\tau=1500$}} \\
\orig & 3.32 & 38 & 2910 & 
$1.03\cdot10^{-6}$ \\
\cmidrule(l){1-5}
\frobbndrk, $C=10$ & 2.37 & 38 & 2487 &
$1.03\cdot10^{-6}$ \\
\frobbndrk, $C=10^{2}$ & 2.36 & 38 & 2462 & 
$1.03\cdot10^{-6}$ \\
\frobbndrk, $C=10^{3}$ & 2.32 & 38 & 2426 &
$1.03\cdot10^{-6}$ \\
\frobbndrk, $C=10^{4}$ & 2.29 & 38 & 2396 & 
$1.03\cdot10^{-6}$ \\
\frobbndrk, $C=10^{5}$ & 2.27 & 38 & 2364 & 
$1.04\cdot10^{-6}$ \\
\frobbndrk, $C=10^{6}$ & 3.64 & 63 & 3846 & 
$4.24\cdot10^{-7}$ \\
\midrule
\midrule
\multicolumn{5}{l}{\emph{Random $5000\times1000$, rank $10$, $60\%$ observed, $\tau=1500$}} \\
\orig & 4.78 & 15 & 608 & 
$9.59\cdot10^{-7}$ \\
\cmidrule(l){1-5}
\frobbndrk, $C=10$ & 2.48 & 15 & 354 & 
$9.59\cdot10^{-7}$ \\
\frobbndrk, $C=10^{2}$ & 2.44 & 15 & 346 & 
$9.59\cdot10^{-7}$ \\
\frobbndrk, $C=10^{3}$ & 2.42 & 15 & 339 & 
$9.59\cdot10^{-7}$ \\
\frobbndrk, $C=10^{4}$ & 2.38 & 15 & 328 & 
$9.72\cdot10^{-7}$ \\
\frobbndrk, $C=10^{5}$ & 2.34 & 15 & 320 & 
$9.73\cdot10^{-7}$ \\
\frobbndrk, $C=10^{6}$ & 3.68 & 26 & 543 & 
$5.60\cdot10^{-7}$ \\
\midrule
\multicolumn{5}{l}{\emph{Random $5000\times1000$, rank $30$, $60\%$ observed, $\tau=1500$}} \\
\orig & 6.94 & 19 & 1316 & 
$9.85\cdot10^{-7}$ \\
\cmidrule(l){1-5}
\frobbndrk, $C=10$ & 6.82 & 19 & 1298 & 
$9.85\cdot10^{-7}$ \\
\frobbndrk, $C=10^{2}$ & 6.73 & 19 & 1280 & 
$9.85\cdot10^{-7}$ \\
\frobbndrk, $C=10^{3}$ & 6.63 & 19 & 1264 & 
$9.85\cdot10^{-7}$ \\
\frobbndrk, $C=10^{4}$ & 6.50 & 19 & 1241 & 
$9.85\cdot10^{-7}$ \\
\frobbndrk, $C=10^{5}$ & 6.66 & 20 & 1271 & 
$5.67\cdot10^{-7}$ \\
\frobbndrk, $C=10^{6}$ & 10.17 & 34 & 2007 & 
$9.15\cdot10^{-7}$ \\
\midrule
\midrule
\multicolumn{5}{l}{\emph{Random $10000\times2000$, rank $10$, $60\%$ observed, $\tau=1500$}} \\
\orig & 18.00 & 9 & 421 & 
$7.88\cdot10^{-7}$ \\
\cmidrule(l){1-5}
\frobbndrk, $C=10$ & 8.02 & 9 & 215 & 
$7.88\cdot10^{-7}$ \\
\frobbndrk, $C=10^{2}$ & 8.02 & 9 & 214 & 
$7.88\cdot10^{-7}$ \\
\frobbndrk, $C=10^{3}$ & 7.88 & 9 & 206 & 
$7.88\cdot10^{-7}$ \\
\frobbndrk, $C=10^{4}$ & 7.80 & 9 & 200 & 
$7.91\cdot10^{-7}$ \\
\frobbndrk, $C=10^{5}$ & 7.59 & 9 & 194 & 
$7.91\cdot10^{-7}$ \\
\frobbndrk, $C=10^{6}$ & 9.78 & 13 & 258 & 
$4.76\cdot10^{-7}$ \\
\midrule
\multicolumn{5}{l}{\emph{Random $10000\times2000$, rank $30$, $60\%$ observed, $\tau=1500$}} \\
\orig & 31.21 & 11 & 1110 & 
$6.66\cdot10^{-7}$ \\
\cmidrule(l){1-5}
\frobbndrk, $C=10$ & 20.58 & 11 & 839 & 
$6.66\cdot10^{-7}$ \\
\frobbndrk, $C=10^{2}$ & 20.30 & 11 & 829 & 
$6.66\cdot10^{-7}$ \\
\frobbndrk, $C=10^{3}$ & 19.97 & 11 & 817 & 
$6.66\cdot10^{-7}$ \\
\frobbndrk, $C=10^{4}$ & 19.70 & 11 & 803 & 
$6.66\cdot10^{-7}$ \\
\frobbndrk, $C=10^{5}$ & 19.17 & 11 & 783 & 
$7.16\cdot10^{-7}$ \\
\frobbndrk, $C=10^{6}$ & 20.36 & 13 & 844 & 
$9.30\cdot10^{-7}$ \\
\bottomrule
\end{tabular}
\caption{SVT equipped with either \orig or \frobbndrk\ on random matrices, for the values of $C$ in the threshold $C/k^{2+\epsilon}$ indicated in the table. }
\label{tab:random_grandi_frobbndrk_C_ranks}
\end{table}

%% file: tabelle_det/tabella_random_grandi_frobbndrk_C_ranghi10e30_obs30_dfac12.tex
\begin{table}[t]
\centering
\small
\begin{tabular}{lcccc}
\toprule
 & Time (s) & SVT it. & $\ell_{\rm tot}$ & 
 $\|X-A\|_F/\|A\|_F$ \\
\midrule
\multicolumn{5}{l}{\emph{Random $2000\times600$, rank $10$, $30\%$ observed, $\tau=2200$}} \\
\orig & 1.80 & 54 & 2478 & 
$1.35\cdot10^{-6}$ \\
\cmidrule(l){1-5}
\frobbndrk, $C=10$ & 0.93 & 54 & 1508 & 
$1.17\cdot10^{-6}$ \\
\frobbndrk, $C=10^{2}$ & 0.88 & 54 & 1474 & 
$1.17\cdot10^{-6}$ \\
\frobbndrk, $C=10^{3}$ & 0.87 & 54 & 1432 & 
$1.17\cdot10^{-6}$ \\
\frobbndrk, $C=10^{4}$ & 0.84 & 54 & 1394 & 
$1.17\cdot10^{-6}$ \\
\frobbndrk, $C=10^{5}$ & 0.81 & 54 & 1356 & 
$1.17\cdot10^{-6}$ \\
\frobbndrk, $C=10^{6}$ & 0.85 & 57 & 1390 & 
$1.06\cdot10^{-6}$ \\
\midrule
\multicolumn{5}{l}{\emph{Random $2000\times600$, rank $30$, $30\%$ observed, $\tau=2200$}} \\
\orig & 4.68 & 94 & 8231 & 
$1.61\cdot10^{-6}$ \\
\cmidrule(l){1-5}
\frobbndrk, $C=10$ & 4.41 & 94 & 8027 & 
$1.59\cdot10^{-6}$ \\
\frobbndrk, $C=10^{2}$ & 4.31 & 94 & 7907 & 
$1.59\cdot10^{-6}$ \\
\frobbndrk, $C=10^{3}$ & 4.15 & 94 & 7796 & 
$1.59\cdot10^{-6}$ \\
\frobbndrk, $C=10^{4}$ & 4.09 & 94 & 7680 & 
$1.59\cdot10^{-6}$ \\
\frobbndrk, $C=10^{5}$ & 3.95 & 94 & 7567 & 
$1.59\cdot10^{-6}$ \\
\frobbndrk, $C=10^{6}$ & 4.07 & 95 & 7521 & 
$1.48\cdot10^{-6}$ \\
\midrule
\midrule
\multicolumn{5}{l}{\emph{Random $5000\times1000$, rank $10$, $30\%$ observed, $\tau=3000$}} \\
\orig & 4.08 & 32 & 1178 & 
$8.84\cdot10^{-7}$ \\
\cmidrule(l){1-5}
\frobbndrk, $C=10$ & 2.95 & 31 & 823 & 
$9.85\cdot10^{-7}$ \\
\frobbndrk, $C=10^{2}$ & 2.91 & 31 & 808 & 
$9.85\cdot10^{-7}$ \\
\frobbndrk, $C=10^{3}$ & 2.97 & 31 & 782 & 
$9.85\cdot10^{-7}$ \\
\frobbndrk, $C=10^{4}$ & 2.92 & 31 & 768 & 
$9.85\cdot10^{-7}$ \\
\frobbndrk, $C=10^{5}$ & 2.75 & 31 & 746 & 
$9.88\cdot10^{-7}$ \\
\frobbndrk, $C=10^{6}$ & 3.35 & 40 & 936 & 
$4.48\cdot10^{-7}$ \\
\midrule
\multicolumn{5}{l}{\emph{Random $5000\times1000$, rank $30$, $30\%$ observed, $\tau=3000$}} \\
\orig & 10.30 & 47 & 3451 & 
$1.23\cdot10^{-6}$ \\
\cmidrule(l){1-5}
\frobbndrk, $C=10$ & 10.06 & 47 & 3452 & 
$1.23\cdot10^{-6}$ \\
\frobbndrk, $C=10^{2}$ & 9.97 & 47 & 3411 & 
$1.23\cdot10^{-6}$ \\
\frobbndrk, $C=10^{3}$ & 9.93 & 47 & 3363 & 
$1.23\cdot10^{-6}$ \\
\frobbndrk, $C=10^{4}$ & 9.80 & 47 & 3321 & 
$1.23\cdot10^{-6}$ \\
\frobbndrk, $C=10^{5}$ & 9.66 & 47 & 3276 & 
$1.23\cdot10^{-6}$ \\
\frobbndrk, $C=10^{6}$ & 9.59 & 48 & 3289 & 
$1.15\cdot10^{-6}$ \\
\midrule
\midrule
\multicolumn{5}{l}{\emph{Random $10000\times2000$, rank $10$, $30\%$ observed, $\tau=4000$}} \\
\orig & 19.92 & 20 & 925 & 
$8.23\cdot10^{-7}$ \\
\cmidrule(l){1-5}
\frobbndrk, $C=10$ & 9.16 & 20 & 496 & 
$7.02\cdot10^{-7}$ \\
\frobbndrk, $C=10^{2}$ & 8.56 & 20 & 485 & 
$7.02\cdot10^{-7}$ \\
\frobbndrk, $C=10^{3}$ & 8.21 & 20 & 475 & 
$7.02\cdot10^{-7}$ \\
\frobbndrk, $C=10^{4}$ & 8.14 & 20 & 464 & 
$7.02\cdot10^{-7}$ \\
\frobbndrk, $C=10^{5}$ & 8.09 & 20 & 453 & 
$7.06\cdot10^{-7}$ \\
\frobbndrk, $C=10^{6}$ & 8.24 & 21 & 462 & 
$5.24\cdot10^{-7}$ \\
\midrule
\multicolumn{5}{l}{\emph{Random $10000\times2000$, rank $30$, $30\%$ observed, $\tau=4000$}} \\
\orig & 29.21 & 26 & 2152 & 
$1.00\cdot10^{-6}$ \\
\cmidrule(l){1-5}
\frobbndrk, $C=10$ & 22.10 & 26 & 1859 & 
$9.23\cdot10^{-7}$ \\
\frobbndrk, $C=10^{2}$ & 22.11 & 26 & 1842 & 
$9.23\cdot10^{-7}$ \\
\frobbndrk, $C=10^{3}$ & 21.83 & 26 & 1814 & 
$9.23\cdot10^{-7}$ \\
\frobbndrk, $C=10^{4}$ & 21.64 & 26 & 1795 & 
$9.23\cdot10^{-7}$ \\
\frobbndrk, $C=10^{5}$ & 21.43 & 26 & 1773 & 
$9.24\cdot10^{-7}$ \\
\frobbndrk, $C=10^{6}$ & 20.92 & 26 & 1740 & 
$1.08\cdot10^{-6}$ \\
\bottomrule
\end{tabular}
\caption{SVT equipped with either \orig or \frobbndrk\ on random matrices, for the values of $C$ in the threshold $C/k^{2+\epsilon}$ indicated in the table. }
\label{tab:random_grandi_frobbndrk_C_ranks_obs30_dfac12}
\end{table}

%% file: Main_paper_deterministic.bbl
\begin{thebibliography}{10}
	\providecommand{\url}[1]{{#1}}
	\providecommand{\urlprefix}{URL }
	\expandafter\ifx\csname urlstyle\endcsname\relax
	\providecommand{\doi}[1]{DOI~\discretionary{}{}{}#1}\else
	\providecommand{\doi}{DOI~\discretionary{}{}{}\begingroup
		\urlstyle{rm}\Url}\fi
	
	\bibitem{jacek1}
	Bellavia, S., Gondzio, J., Porcelli, M.: An inexact dual logarithmic barrier
	method for solving sparse semidefinite programs.
	\newblock Mathematical Programming \textbf{178}(1), 109--143 (2019)
	
	\bibitem{jacek2}
	Bellavia, S., Gondzio, J., Porcelli, M.: A relaxed interior point method for
	low-rank semidefinite programming problems with applications to matrix
	completion.
	\newblock Journal of Scientific Computing \textbf{89}(2), 46 (2021)
	
	\bibitem{bellavia2025inexact}
	Bellavia, S., Rebegoldi, S., Silei, M.: An inexact alternating projection
	method with application to matrix completion.
	\newblock Computational Optimization and Applications \textbf{95}, 29--74
	(2026)
	
	\bibitem{boumal2015low}
	Boumal, N., Absil, P.A.: Low-rank matrix completion via preconditioned
	optimization on the grassmann manifold.
	\newblock Linear Algebra and its Applications \textbf{475}, 200--239 (2015)
	
	\bibitem{cai2010singular}
	Cai, J.F., Cand{\`e}s, E.J., Shen, Z.: A singular value thresholding algorithm
	for matrix completion.
	\newblock SIAM Journal on optimization \textbf{20}(4), 1956--1982 (2010)
	
	\bibitem{CandesRecht2009}
	Cand{\`e}s, E.J., Recht, B.: Exact matrix completion via convex optimization.
	\newblock Foundations of Computational Mathematics \textbf{9}(6), 717--772
	(2009)
	
	\bibitem{survey1}
	Chen, Z., Wang, S.: A review on matrix completion for recommender systems.
	\newblock Knowledge and Information Systems \textbf{64}(1), 1--34 (2022)
	
	\bibitem{Relu26}
	Gillis, N., Porcelli, M., Seraghiti, G.: An extrapolated and provably
	convergent algorithm for nonlinear matrix decomposition with the relu
	function.
	\newblock SIAM Journal on Mathematics of Data Science \textbf{8}(3), 738--765
	(2026).
	\newblock \doi{10.1137/25M1746859}
	
	\bibitem{GolubKahan65}
	Golub, G., Kahan, W.: {Calculating the Singular Values and Pseudo-Inverse of a
		Matrix}.
	\newblock Journal of the Society for Industrial and Applied Mathematics Series
	B Numerical Analysis \textbf{2}(2), 205--224 (1965).
	\newblock \doi{10.1137/0702016}
	
	\bibitem{golub2013matrix}
	Golub, G.H., Van~Loan, C.F.: Matrix computations.
	\newblock JHU press (2013)
	
	\bibitem{goyens2023nonlinear}
	Goyens, F., Cartis, C., Eftekhari, A.: Nonlinear matrix recovery using
	optimization on the grassmann manifold.
	\newblock Applied and Computational Harmonic Analysis \textbf{62}, 498--542
	(2023)
	
	\bibitem{loraine}
	Habibi, S., Ko{\v{c}}vara, M., Stingl, M.: Loraine--an interior-point solver
	for low-rank semidefinite programming.
	\newblock Optimization Methods and Software \textbf{39}(6), 1185--1215 (2024)
	
	\bibitem{halko2011finding}
	Halko, N., Martinsson, P.G., Tropp, J.A.: Finding structure with randomness:
	Probabilistic algorithms for constructing approximate matrix decompositions.
	\newblock SIAM review \textbf{53}(2), 217--288 (2011)
	
	\bibitem{higham_power_2024}
	Higham, N.: {The power of bidiagonal matrices}.
	\newblock Electron. J. Linear Algebra \textbf{40}, 453–474 (2024).
	\newblock \doi{10.13001/ela.2024.8297}
	
	\bibitem{horn1991topics}
	Horn, R.A., Johnson, C.R.: Topics in matrix analysis  (1991)
	
	\bibitem{kim2020deep}
	Kim, S., Nguyen, L.T., Shim, B.: Deep neural network based matrix completion
	for internet of things network localization.
	\newblock In: ICASSP 2020-2020 IEEE International Conference on Acoustics,
	Speech and Signal Processing (ICASSP), pp. 3427--3431. IEEE (2020)
	
	\bibitem{propack_doc}
	Larsen, R.M.: Lanczos bidiagonalization with partial reorthogonalization.
	\newblock Technical Report DAIMI PB-357, Department of Computer Science, Aarhus
	University (1998).
	\newblock \urlprefix\url{https://tidsskrift.dk/daimipb/article/view/7070}
	
	\bibitem{larsen2004propack}
	Larsen, R.M.: {PROPACK}: Software for large and sparse {SVD} calculations.
	\newblock \url{http://sun.stanford.edu/~rmunk/PROPACK/} (2004)
	
	\bibitem{li_convergence_2015}
	Li, R.C., Zhang, L.H.: Convergence of the block {Lanczos} method for eigenvalue
	clusters.
	\newblock Numerische Mathematik \textbf{131}(1), 83--113 (2015).
	\newblock \doi{10.1007/s00211-014-0681-6}.
	\newblock \urlprefix\url{http://link.springer.com/10.1007/s00211-014-0681-6}
	
	\bibitem{liu2010interior}
	Liu, Z., Vandenberghe, L.: Interior-{Point} {Method} for {Nuclear} {Norm}
	{Approximation} with {Application} to {System} {Identification}.
	\newblock SIAM Journal on Matrix Analysis and Applications \textbf{31}(3),
	1235--1256 (2010).
	\newblock \doi{10.1137/090755436}.
	\newblock \urlprefix\url{http://epubs.siam.org/doi/10.1137/090755436}
	
	\bibitem{fpca}
	Ma, S., Goldfarb, D., Chen, L.: Fixed point and bregman iterative methods for
	matrix rank minimization.
	\newblock Mathematical Programming \textbf{128}(1), 321--353 (2011)
	
	\bibitem{oh2018fast}
	Oh, T.H., Matsushita, Y., Tai, Y.W., Kweon, I.S.: Fast randomized singular
	value thresholding for low-rank optimization.
	\newblock IEEE Transactions on Pattern Analysis and Machine Intelligence
	\textbf{40}(2), 376--391 (2018).
	\newblock \doi{10.1109/TPAMI.2017.2677440}
	
	\bibitem{survey2}
	Ramlatchan, A., Yang, M., Liu, Q., Li, M., Wang, J., Li, Y.: A survey of matrix
	completion methods for recommendation systems.
	\newblock Big Data Mining and Analytics \textbf{1}(4), 308--323 (2018)
	
	\bibitem{recht2010guaranteed}
	Recht, B., Fazel, M., Parrilo, P.A.: Guaranteed {Minimum}-{Rank} {Solutions} of
	{Linear} {Matrix} {Equations} via {Nuclear} {Norm} {Minimization}.
	\newblock SIAM Review \textbf{52}(3), 471--501 (2010).
	\newblock \doi{10.1137/070697835}.
	\newblock \urlprefix\url{http://epubs.siam.org/doi/10.1137/070697835}
	
	\bibitem{saad_rates_1980}
	Saad, Y.: On the {Rates} of {Convergence} of the {Lanczos} and the
	{Block}-{Lanczos} {Methods}.
	\newblock SIAM Journal on Numerical Analysis \textbf{17}(5), 687--706 (1980).
	\newblock \doi{10.1137/0717059}.
	\newblock \urlprefix\url{http://epubs.siam.org/doi/10.1137/0717059}
	
	\bibitem{saad2011numerical}
	Saad, Y.: Numerical methods for large eigenvalue problems: revised edition.
	\newblock SIAM (2011)
	
	\bibitem{silei2023recovering}
	Silei, M., Bellavia, S., Superchi, F., Bianchini, A.: Recovering corrupted data
	in wind farm measurements: A matrix completion approach.
	\newblock Energies \textbf{16}(4), 1674 (2023)
	
	\bibitem{apgl}
	Toh, K.C., Yun, S.: An accelerated proximal gradient algorithm for nuclear norm
	regularized linear least squares problems.
	\newblock Pacific Journal of optimization \textbf{6}(615-640), 15 (2010)
	
	\bibitem{vandereycken2013low}
	Vandereycken, B.: Low-rank matrix completion by riemannian optimization.
	\newblock SIAM Journal on Optimization \textbf{23}(2), 1214--1236 (2013)
	
	\bibitem{wang2021fast}
	Wang, H., Cai, J., Wang, T., Wei, K.: Fast {C}adzow’s algorithm and a
	gradient variant.
	\newblock J. Sci. Comput. \textbf{88}(2), 41 (2021)
	
	\bibitem{wang2025solving}
	Wang, J., Hu, L.: Solving low-rank semidefinite programs via manifold
	optimization.
	\newblock Journal of Scientific Computing \textbf{104}(1), 33 (2025)
	
	\bibitem{LMaFit}
	Wen, Z., Yin, W., Zhang, Y.: Solving a low-rank factorization model for matrix
	completion by a nonlinear successive over-relaxation algorithm.
	\newblock Mathematical Programming Computation \textbf{4}(4), 333--361 (2012)
	
\end{thebibliography}
